\documentclass[a4paper]{article}

\title{Direct embeddings of relatively hyperbolic groups with optimal $\ell^p$ compression exponent}
\author{David Hume \\ \small{University of Oxford} \\ \small{hume@maths.ox.ac.uk}}
\date{}

\newcommand{\symd}{\bigtriangleup}
\newcommand{\set}[1]{\left\{#1\right\}}
\newcommand{\setcon}[2]{\left\{#1\ \left|\ #2\right.\right\}}
\newcommand{\norm}[1]{\left\lVert#1\right\rVert}
\newcommand{\abs}[1]{\left\lvert#1\right\rvert}
\newcommand{\into}{\hookrightarrow}
\newcommand{\R}{\mathbb{R}}
\newcommand{\Z}{\mathbb{Z}}
\newcommand{\N}{\mathbb{N}}
\newcommand{\m}{\medskip \\}
\newcommand{\h}{\hspace{3mm}}
\newcommand{\tu}{\textup}
\newcommand{\pf}{\noindent {\bf Proof: \ \ }}

\newcommand{\geo}[1]{\underline{#1}}

\newcommand{\setg}[2]{[\![#1,#2]\!]}

\usepackage{amssymb}
\usepackage{amsmath, amssymb, amsthm}
\usepackage{graphicx}
\usepackage{enumitem}
\usepackage{mhsetup, mathtools}
\usepackage[utf8]{inputenc}
\usepackage{tikz}
\usepackage{MnSymbol}
\usepackage{float}
\usepackage{setspace}

\newtheorem{bthm}{Theorem}
\newtheorem{thm}{Theorem}[section]
\newtheorem{cor}[thm]{Corollary}
\newtheorem{lem}[thm]{Lemma}
\newtheorem{defn}[thm]{Definition}
\newtheorem{prop}[thm]{Proposition}

\begin{document}
\maketitle
\begin{abstract}
 \noindent We prove that for all $p>1$, every relatively hyperbolic group has $\ell^p$ compression exponent equal to the minimum of the exponents of its maximal peripheral subgroups. This improves results of Dadarlat-Guentner and Dreesen. As a first step we give a direct geometric proof that hyperbolic groups have $\ell^p$ compression exponent $1$, independent of those given by Bonk-Schramm, Buyalo-Dranishnikov-Schroeder, Gal and Tessera.
\end{abstract}

\section{Introduction and statement of results}
Coarse embeddings of discrete metric spaces into Banach spaces have been an important topic in computer science for many years \cite{LLR95,HLW06,BDGRRRS} and more recently became so in geometric group theory, combinatorics and $K$-theory. At the suggestion of Gromov, Yu and later Kasparov and Yu proved that any group admitting a coarse embedding into any uniformly convex Banach space satisfies the Novikov and coarse Baum-Connes conjectures \cite{Gr93,Yu00,KY06}. However, Gromov proved that there exist finitely generated groups which do not admit a coarse embedding into any Hilbert space \cite{Gr00,AD08}.
\m
We recall that a \emph{coarse \tu{(}or uniform\tu{)} embedding} of $(X,d)$ into $(Y,d')$ is a map $\phi:X\to Y$ for which there are functions $\rho_\pm:\mathbb{R}_{\geq 0}\to\mathbb{R}_{\geq 0}$, such that $\rho_\pm(r)\to\infty$ as $r\to\infty$ and 
$\rho_-\left(d(x_1,x_2)\right)\leq d'(\phi(x_1),\phi(x_2)) \leq \rho_+\!\left(d(x_1,x_2)\right)$.
\\
A stronger notion was introduced by Guentner and Kaminker, called the compression exponent, $\alpha^*_Y(X)$ \cite{GK04}. This is defined as the supremum over $\alpha\in[0,1]$ such that there exists a coarse embedding $\phi:X\to Y$ with $\rho_+(r)\leq Cr+C$ and $\rho_-(r)\geq C^{-1}r^\alpha-C$. We denote by $\alpha^*_p(X)$ the value $\alpha^*_{\ell^p(\N)}(X)$, and call this value the $\ell^p$ compression exponent of $X$.
\m
The same definitions can be made for equivariant embeddings, that is embeddings $\phi:G\to X$ such that for some action of $G$ on $X$, $(g,x) \mapsto g\cdot x$, $\phi$ is the orbit map of the action with respect to some point $e\in X$, i.e. $\phi(g)=g\cdot e$ for every $g\in G$. The equivariant $\ell^p$ compression exponent $\alpha^\#_p(G)$ is then given by the same supremum as $\alpha^*_p(G)$ but taken over only equivariant embeddings. While this value is independent of the choice of finite generating set for a group, it is not known to be a quasi-isometry invariant (unlike the non-equivariant compression exponent which is obviously such an invariant).
\m
The $\ell^p$ compression exponents of a finitely generated group, $\alpha_p^*(G)$ and $\alpha^\#_p(G)$ are closely linked to various forms of amenability, the speed of random walks on the Cayley graph and notions of large-scale (asymptotic) dimension.
\\
We recall that the speed of random walks in a group $G$, $\beta^*(G)$, is the supremum over $\beta\in\left[0,1\right]$ for which there exists a constant $C$ such that the expected length of any concatenation of $n$ generators of $G$ is at least $C^{-1}n^\beta - C$. This is directly linked to $\ell^p$ compression via the inequalities $\alpha^*_p(G)\leq \frac{1}{q\beta^*(G)}$ which holds for all amenable groups and $\alpha^\#_p(G)\leq \frac{1}{q\beta^*(G)}$, which holds in general \cite{NP08}. In these inequalities it is important to note that the value $q=\min\set{p,2}$ is the power type of the modulus of smoothness, which for $\ell^p$ spaces with $p>1$ is equal to $\min\set{p,2}$.
\\
In particular, this result says that any group satisfying $\alpha_p^\#(G)>\frac{1}{2}$ for some $p\geq 2$ is amenable as every non-amenable group has $\beta^*(G)=1$. However, there are solvable groups with $\alpha^\#_p(G)=0$ for all $p$ \cite{GK04,NP08,Au11}. Groups with $\alpha_p^*(G)>\frac{1}{p}$ for some $p\in [1,2]$ satisfy property (A), a non-equivariant form of amenability \cite{Yu00,GK04,WW75}.
\m
In \cite{Gr87}, Gromov introduced relatively hyperbolic groups as a generalisation of hyperbolic groups. The class of relatively hyperbolic groups includes: hyperbolic groups, amalgamated products and HNN-extensions over finite subgroups, fully residually free (limit) groups \cite{Da03,Al05} - which are key objects in solving the Tarski conjecture \cite{Se01,KM10}, geometrically finite Kleinian groups and fundamental groups of non-geometric closed $3$-manifolds with at least one hyperbolic component \cite{Da03}. 
\\
These groups have many different characterisations: in terms of group actions \cite{Bo99}, group-theoretic structure \cite{Fa98}, \cite{Da03'}, \cite{Os06}, dynamics on the boundary \cite{Ya04} and metric geometry \cite{DS05}. 
\m
Our main result is the following:
\begin{bthm}\label{RHgps} \tu{(cf. $\ref{RHgps'}$)} \h  Let $G$ be a finitely generated group which is hyperbolic relative to a collection of finitely generated subgroups $\setcon{H_i}{i=1,\dots,n}$. For all $p>1$, $\alpha^*_p(G)=\min\setcon{\alpha_p^*(H_i)}{i=1,\dots,n}$.
\end{bthm}

\noindent In fact, we obtain this result in greater generality (Theorem $\ref{atg'}$), and use it to show $\alpha_p^*(\pi_1(M))=1$ for any closed $3$-manifold $M$ and any $p>1$ (Corollary $\ref{3mgps}$). We also obtain bounds on the more commonly studied $L_p$ compression (Corollary $\ref{L_p}$).
\m
The lower bound $\alpha^*_p(G)\leq\min\set{\alpha_p^*(H_i)}$ in Theorem $\ref{RHgps}$ is clear as the subgroups $H_i$ are quasi-isometrically embedded in $G$.
\m
Certain specific subclasses of relatively hyperbolic groups were previously known satisfy the equality in Theorem \ref{RHgps}. Sela proved the coarse embeddability of hyperbolic groups into Hilbert spaces \cite{Se96}, but his direct methods only yield a lower bound of $\frac{1}{2}$ on compression exponent. For free groups, compression exponent $1$ was determined by Guentner and Kaminker, with refinements by Brodskiy and Higes and optimal results obtained independently by Gal and Tessera \cite{GK04,BS08,Gal,Te11}. A lower bound of $1$ for general hyperbolic groups can be obtained using major theorems from \cite{BoS} and \cite{BS05} or from \cite{Te11}.
\m
All relatively hyperbolic groups admit coarse embeddings into Hilbert spaces (providing their maximal peripheral subgroups do) \cite{DG07}, (previous results \cite{CDGY,DG03}), but the methods here do not generalise to other $\ell^p$ spaces and cannot be used to provide a lower bound on compression exponent so a completely different technique is required to progress. However, bounds are available in restricted cases. Dreesen proves that given finitely generated groups $A$, $B$ and $C$, where $C$ is a finite subgroup of $A$ and $B$, $\min\left\{\alpha_2^*(A),\alpha_2^*(B),\frac{1}{2}\right\}\leq\alpha_2^*(A*_CB)$ and $\min\left\{\alpha_2^*(A),\frac{1}{2}\right\}\leq\alpha_2^*(\tu{HNN}(A,C,\theta))$ \cite{Dr11}. 
\m
For limit groups the result is also already known. It follows from work of Wise that limit groups quasi-isometrically embed into right-angled Artin groups, so these specific relatively hyperbolic groups have Hilbert compression exponent $1$ \cite{Wi11,DJ00,DJ99,GK04}.
\m
{\bf Outline:} \h Theorem \ref{RHgps} is technically difficult so we approach it via new results for two key sub-collections, the methods presented aid the intuition and strategy behind the final proof. In Section \ref{hyp1} we provide a new, direct, self-contained proof that the $\ell^p$ compression of any hyperbolic group is $1$. This is not a new result, but all previous proofs rely in a key way on major theorems of Bonk-Schramm or Buyalo-Schroeder \cite{BoS, BS05}. It is the fact that this embedding is direct and easily constructed which is the main interest here as it allows the possibility of extending such a result to relatively hyperbolic spaces. Then we find embeddings of amalgamated products and HNN-extensions over finite groups displaying the exact compression exponent (Section \ref{Freepf}), by applying a careful weighting procedure to Dreesen's method. 
A na\"ive implementation of these two methods will never yield a Lipschitz embedding in the general case, so for the final proof (Section \ref{mainthm}) we are required to refine these two arguments and ensure that a suitable combination of them also provides an optimal lower bound on compression exponent.
\m
{\bf Acknowledgements:}\h The author would like to thank Romain Tessera and an anonymous referee for comments on previous versions of this paper and is grateful for the support of the EPSRC through a D.Phil.\! student grant and the grant ``Geometric and analytic aspects of infinite groups''.

\section{Preliminaries}
Firstly, we will use the following notation throughout the paper.
\m
{\bf Notation:}\h We will define $\setg{x}{y}$ to be the set of all geodesics from $x$ to $y$ in a given metric space.
\m
Given functions $f,g:\N\to\R$ we write $f\preceq g$ to mean that there exists some constant $C$ such that $f(x)\leq Cg(x)+C$ for all $x\in\R$. We write $f\asymp g$ if $f\preceq g$ and $g\preceq f$.
\m
In this section we define the collection of functions that will be used throughout this paper as lower bounds to the coarse embeddings we construct and present two lemmas which will be used in later sections.

\begin{defn}\label{defnccp} \h {\bf Concave functions and property} \tu{($C^c_p$)} 
\m
We will call a function $f:\N\to\R_{\geq 0}$ concave if $f$ is non-decreasing and for all $m,n\in\N$ with $n\geq m$:
\[
f(n+m)-f(n)\leq f(n)-f(n-m).
\]
Let $f:\N\to\R_{\geq 0}$ be a concave function satisfying Tessera's property \tu{($C_p$)},
\[
  \sum_{n=1}^\infty \frac{1}{n}\left(\frac{f(n)}{n}\right)^p < \infty.
\]
$f$ is said to satisfy \tu{($C^c_p$)} if, in addition, $\frac{f(n)^p}{n}$ is non-decreasing for all $n$ sufficiently large.
\end{defn}
\noindent We observe here that for all $\epsilon>0$ and all $p>1$,
\[
f(n)=\frac{n}{(\log_2(n+2)(\log_2\log_2)^{1+\epsilon}(n+2))^\frac{1}{p}} 
\]
has property ($C^c_p$).
\m
The concavity condition above is modelled on the usual concavity condition $f''\leq 0$ given for smooth functions.
\m
We now present two technical lemmas which will later explain the necessity of the preceding definition.

\begin{lem}\label{lbound} \h Let $M$ be a finite subset of $\N$ such that $M=\{m_1,m_2,\dots, m_{2k}\}$ with $m_i < m_{i+1}$ and $m_1\geq 1$.
\m
Let $p>1$ and let $f:\N\to\R_{\geq 0}$ be a concave function such that $\frac{f(n)^p}{n}$ is non-decreasing. Then
\[
 \sum_i \frac{f(m_{2i})^p}{m_{2i}}(m_{2i}-m_{2i-1}) \geq \left(\frac{1}{2}\right)^{3+p} f\left(\sum_i m_{2i}-m_{2i-1}\right)^p.
\]
\end{lem}
\pf For ease of notation we set $\displaystyle m= \sum_{i=1}^k m_{2i}-m_{2i-1}$.
\m
As $\frac{f(n)^p}{n}$ is non-decreasing,
\[
\sum_{i=1}^k \frac{f(m_{2i})^p}{m_{2i}}(m_{2i}-m_{2i-1}) \geq \sum_{n=1}^{m} \frac{f(n)^p}{n}.
\]
The result then follows from the method in \cite[Theorem $7.3$]{Te11}.
\[
\begin{array}{rcl}
\displaystyle \sum_{n=1}^{m} \frac{f(n)^p}{n}
 &
\geq
 &
\displaystyle\sum_{n=m/2}^{m} \frac{1}{n} f([m/2])^p
\m
 &
\geq 
 &
\displaystyle\frac{1}{4} f([m/2])^p \geq \frac{1}{2}^{3+p} f(m).
\end{array}
\]

\begin{lem}\label{ubound} \h Let $M$ be a finite subset of $\N$ such that $M=\{m_1,m_2,\dots, m_{2k}\}$ with $m_i < m_{i+1}$ and $m_1\geq 1$.
\m
Let $p>1$ and let $f:\N\to\R_{\geq 0}$ be a concave function with property \tu{($C_p$)}. Then there exists some uniform constant $C$ such that
\[
\sum_i \left(\frac{f(m_{2i})}{m_{2i}}\right)^p\frac{m_{2i}-m_{2i-1}}{m_{2i}} \leq C.
\]
Moreover, if for each $i$, $m_{2i}\leq 2m_{2i-1}$, then 
 \[
\sum_i \left(\frac{f(m_{2i-1})}{m_{2i-1}}\right)^p\frac{m_{2i}-m_{2i-1}}{m_{2i-1}} \leq 2^{p+1}C.
\]
\end{lem}
\pf As $f$ is concave $\frac{f(n)}{n}$ is non-increasing. Hence 
\[
 \sum_i \left(\frac{f(m_{2i})}{m_{2i}}\right)^p (m_{2i}-m_{2i-1}) \leq \sum_{n=m_{2i-1}+1}^{m_{2i}} \left(\frac{f(n)}{n}\right)^p.
\]
Therefore,
\[
 \sum_i \left(\frac{f(m_{2i})}{m_{2i}}\right)^p\frac{m_{2i}-m_{2i-1}}{m_{2i}} \leq \sum_{n=1}^\infty \frac{1}{n} \left(\frac{f(n)}{n}\right)^p
\]
which is uniformly bounded as $f$ has property \tu{($C_p$)}.
\m
For the second part just notice that $\frac{m_{2i}-m_{2i-1}}{m_{2i-1}} \leq 2\frac{m_{2i}-m_{2i-1}}{m_{2i}}$ and as $f$ is non-decreasing,
\[
 \left(\frac{f(m_{2i-1})}{m_{2i-1}}\right)^p \leq 2^p\left(\frac{f(m_{2i})}{m_{2i}}\right)^p.
\]
\qed

\section{Hyperbolic metric spaces}\label{hyp1}
In this section we provide a short, self-contained and explicit method of embedding uniformly discrete hyperbolic metric spaces with bounded geometry into $\ell^p$ spaces with optimal compression exponent. This builds on ideas used by Tu to prove hyperbolic groups have property (A) \cite{Tu01}.
\m
We first require the following basic lemma in hyperbolic geometry.

\begin{lem}\label{hypstab} \h  Let $X$ be a $\delta$-hyperbolic metric space and let $e\in X$. Let $n\geq 3\delta$ and let $x,y\in X$ with $d(x,e)\geq n$ and $d(x,y)\leq\frac{n}{4}$. For all geodesics $\geo{g_0}\in\setg{x}{e}$, $\geo{g}\in\setg{y}{e}$ and points $p\in\geo{g}([n,2n])$,
\[
 d\left(p,\geo{g_0}([n/2, 5n/2])\right)\leq 3\delta.
\]
\end{lem}
\pf We use the Rips definition of hyperbolicity, so in a geodesic triangle any edge is contained in the union of the $\delta$-neighbourhoods of the other two. Select $p\in\geo{g}([n,2n])$, if $p$ lies within the $\delta$-neighbourhood of $\geo{g_0}$ then we are done as a sufficiently close point must lie within the required range. 
\m
Alternatively, $p$ must lie within the $\delta$-neighbourhood of any geodesic in $\setg{x}{y}$.
\m
Let $z$ be a point on some geodesic in $\setg{x}{y}$ with $d(p,z)<\delta$, then 
\[
d(x,p)\leq d(x,z)+d(z,p)\leq \frac{n}{4} + \delta.
\]
However, $d(x,p)\geq d(y,p)-d(x,y) \geq \frac{3n}{4}$, which is a contradiction as $n\geq 3\delta$. \qed
\begin{figure}[H]
 \centering
  \begin{tikzpicture}[yscale=2.7,xscale=1.8, 
vertex/.style={draw,fill,circle,inner sep=0.3mm}]

\filldraw[black!20!white, draw=black!22!white, fill opacity=0.7]
(-1, 0.2) arc(90:-90:0.2cm) -- (-5,-0.2) arc(270:90:0.2cm) -- (-1,0.2);

\node[vertex]
(c) at ( -6, 0) {};

\filldraw[draw=black!50!white, fill=black!35!white, fill opacity=0.7] (0,0) circle (5mm);

\filldraw%
[%
draw=black!35!white,
very thin,%
fill=black!15!white,%
fill opacity=0.6
]
(120:5 mm)
.. controls +(210:6.5 mm) and (-1.8+0.5,0) ..
(c.east)
.. controls (-1.8+0.5,0) and +(150:6.5 mm) ..
(240:5 mm)
arc (240:120:5 mm);

\begin{scope}
\clip (-1.5, 0.2) arc(90:-90:0.2cm) -- (-4.5,-0.2) arc(270:90:0.2cm) -- (-1.5,0.2);
\filldraw%
[%
draw=black!35!white,
very thin,%
fill=black!65!white,%
fill opacity=0.6
]
(120:5 mm)
.. controls +(210:6.5 mm) and (-1.8+0.5,0) ..
(c.east)
.. controls (-1.8+0.5,0) and +(150:6.5 mm) ..
(240:5 mm)
arc (240:120:5 mm);
\end{scope}

\node[vertex]
(c) at ( -6, 0) {};
\path (c) node[below left] {$e$};

\node[vertex]
(x) at ( 0, 0) {};
\path (x) node[below] {$x$};

\path (0,0.03) node[above] {$B(x;\frac{n}{4})$};

\draw[black] (c) -- (x);
\path (-0.65,0) node[below] {$\geo{g_0}$};

\draw[black!50!white]
(90:4.5 mm)
.. controls +(220:6.5 mm) and (-1.8+0.5,0) ..
(c.east);
\begin{scope}
\clip (-1.5, 0.2) arc(90:-90:0.2cm) -- (-4.5,-0.2) arc(270:90:0.2cm) -- (-1.5,0.2);
\draw[black]
(90:4.5 mm)
.. controls +(220:6.5 mm) and (-1.8+0.5,0) ..
(c.east);
\end{scope}
\node[vertex]
(y) at ( 90:4.5mm) {};
\path (y) node[below right] {$y$};

\draw[black, <->, thin] (0,-0.24) -- (-2,-0.24);
\path (-1,-0.24) node[below] {$n$};
\draw[black, <->, thin] (-2,-0.24) -- (-4,-0.24);
\path (-3,-0.24) node[below] {$n$};

\path (-4,0.1) node[] {$N(\geo{g_0}([n/2,5n/2]);3\delta)$};
\end{tikzpicture}
 \caption{The conclusion of Lemma \ref{hypstab}}\label{fighypstab}
\end{figure}
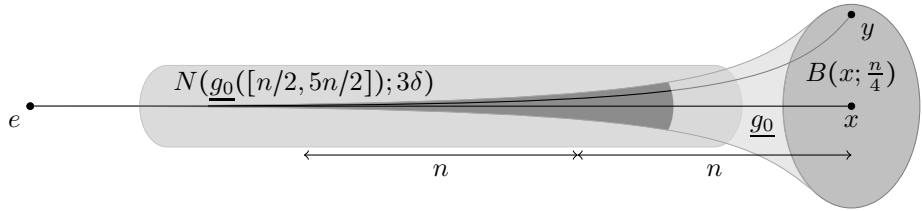

\begin{thm} \label{hyp'} \h Let $X$ be a countable uniformly discrete $\delta$-hyperbolic metric space with bounded geometry. Then given any $p\geq 1$ and any concave function $f$ with property \tu{($C_p$)} there exists a map $\phi:X\to\bigoplus_n\ell^p(X)$ such that for all $x,y\in X$,
\[
 f(d_X(x,y)) \preceq \norm{\phi(x)-\phi(y)}_p \preceq d_X(x,y).
\]
In particular, for all $p\geq 1$, $\alpha_p^*(X)=1$.
\end{thm}
\noindent Throughout this paper, any direct sum of $\ell^p$ spaces is equipped with the $\ell^p$ norm, so all such spaces are isometric to $\ell^p(\N)$. We include the additional detail to more clearly define how each embedding is constructed.
\m
\pf We can reduce our problem to the case where $X$ is the $0$-skeleton of a connected simplicial graph, using \cite[Lemmas $4.1$ and $7.3$]{Tu01}. As $X$ has bounded geometry we can define $N(k)$ to be a bound on the cardinality of any ball of radius $k$.
\m
Fix a basepoint $e\in X$. Given the following collection of restricted geodesics
\[
\geo{G_{x,k,n}}\ \coloneqq \setcon{\geo{g}([n,2n])}{\geo{g}\in\setg{y}{e}\ \tu{for some}\ y\ \tu{with}\ d(x,y)\leq k},
\]
we define $F_{x,k,n}$ to be the set of all points in $X$ lying on some $\geo{g}\in\geo{G_{x,k,n}}$ but not in $B(e;3\delta)$ and set $F(x,k,n)$ to be the characteristic function of $F_{x,k,n}$. $F_{x,k,n}$ is a finite set, so $F(x,k,n)\in\ell^p(X)$.
\m
We then average these functions over all suitable values of $k$, 
\[
H(x,n)=\frac{1}{n}\sum_{k\leq\frac{n}{4}} F(x,k,n).
\]
For increasing values of $k$, the collection of all points lying on a geodesic in the set $\geo{G_{x,k}}\coloneqq\setcon{\geo{g}\in\setg{y}{e}}{d_X(x,y)\leq k}$ of geodesics from the ball of radius $k$ around $x$ to $e$ forms a sequence of `trumpets'.
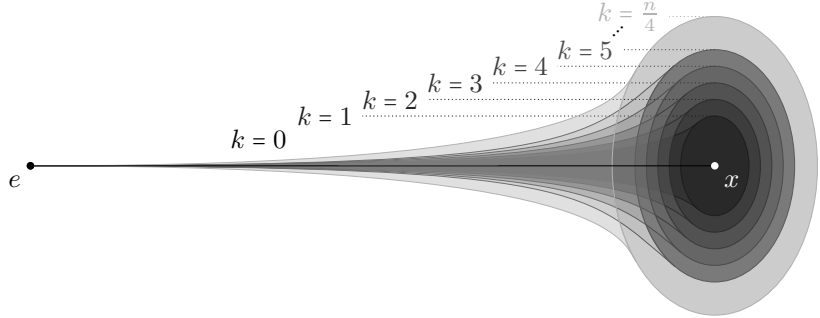
\begin{figure}[H]
\centering
\begin{tikzpicture}[yscale=2.2,xscale=1.5, 
vertex/.style={draw,fill,circle,inner sep=0.3mm}]
\fill[white, fill opacity=0] (-6.22,-0.95) rectangle (0.9,1);
\node[vertex]
(c) at ( -6, 0) {};
{
\filldraw[draw=black!30,fill=black!20] (0,0) circle (0.9cm);
}
{
\node at (-0.85,0.81) {\scriptsize$\udots$};
}
{
\foreach \r / \col / \labeltext in%
{
7/black!40!white/{$k=5$},%
6/black!50!white/{$k=4$},%
5/black!60!white/{$k=3$},%
4/black!70!white/{$k=2$},%
3/black!80!white/{$k=1$}%
}
{%
{%
\begin{scope}
\filldraw%
[%
draw=\col!60!black,
fill=\col!80!black
]
(0,0) circle (\r*0.1 cm);
\end{scope}
\draw[\col!50!black,densely dotted]
(0,\r*0.1 cm) -- +(-4.8 + 0.57*\r,0)
node[pos=1,left] {\labeltext};
}
}
{%
\draw[black!30!white,densely dotted]
(0,0.9 cm) -- +(-0.4,0)
node[pos=1,left] {$k=\frac{n}{4}$};
}
\foreach \r / \col in%
{3/black!80!white, 4/black!70!white, 5/black!60!white, 6/black!50!white, 7/black!40!white}{%
{
\filldraw%
[%
draw=\col!60!black,
very thin,%
fill=\col,%
fill opacity=1.6-0.2*\r
]
(120:\r mm)
.. controls +(210:1.3*\r mm) and (-1.8+\r*0.1,0) ..
(c.east)
.. controls (-1.8+\r*0.1,0) and +(150:1.3*\r mm) ..
(240:\r mm)
arc (240:120:\r mm);
}
}
}
{
\filldraw[draw=black!30!white,very thin,fill=black!40,fill opacity=0.3]%
(140:9mm)
.. controls  +(230:2mm) and (-1,0) ..
(c.east)
.. controls (-1,0) and +(130:2mm) ..
(220:9mm)
arc (220:140:9mm);
}
{
\draw[black]
(c) -- (0,0);
}
{
\draw[black!70!white, very thin]
(-3.3,0) -- (-1.9,0);
}
\node[vertex]
(c) at ( -6, 0) {};
\path (c) node[below left] {$e$};
{
\node[vertex, color=white]
(X) at ( 0, 0) {};
\path (X) node[below right, color=white] {$x$};
}
{
\path ( -4, 0.05) node[above] {$k=0$};
}
\end{tikzpicture}
\caption{A weighted sum of hyperbolic `trumpets'}\label{fighyptrump}
\end{figure}
\noindent The restriction of these to the desired interval (after rescaling) gives the function $H(x,n)$:
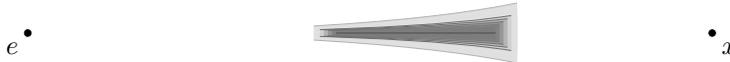
\begin{figure}[H]
\centering
\begin{tikzpicture}[yscale=2,xscale=1.5, 
vertex/.style={draw,fill,circle,inner sep=0.3mm}]
\clip (-6.22,-0.2) rectangle (0.9,0.2);
{
\fill[white, fill opacity=0] (-6.22,-1) rectangle (0.9,1);
\node[vertex]
(c) at ( -6, 0) {};

{
\filldraw[draw=black!30,fill=black!20] (0,0) circle (0.9cm);
}
{
\foreach \r / \col / \labeltext in%
{
7/black!40!white/{$k=5$},%
6/black!50!white/{$k=4$},%
5/black!60!white/{$k=3$},%
4/black!70!white/{$k=2$},%
3/black!80!white/{$k=1$}%
}{
\begin{scope}
\filldraw%
[%
draw=\col!60!black,
fill=\col!80!black
]
(0,0) circle (\r*0.1 cm);
\end{scope}
}

]
\foreach \r / \col in%
{3/black!80!white, 4/black!70!white, 5/black!60!white, 6/black!50!white, 7/black!40!white}{%
{
\filldraw%
[%
draw=\col!60!black,
very thin,%
fill=\col,%
fill opacity=1.6-0.2*\r
]
(120:\r mm)
.. controls +(210:1.3*\r mm) and (-1.8+\r*0.1,0) ..
(c.east)
.. controls (-1.8+\r*0.1,0) and +(150:1.3*\r mm) ..
(240:\r mm)
arc (240:120:\r mm);
}
{
\filldraw%
[%
draw=white,
very thin,%
fill=white,%
]
(200:\r + 32mm) arc(200:160:\r + 32mm)
 -- (-4,1) 
 -- (-6,1) -- (-6,-1) -- 
 (-4,-1)
 -- cycle;
\filldraw%
[%
draw=white,
very thin,%
fill=white,%
]
(c) +(20:\r + 40mm) arc(20:-20:\r + 40mm)
 -- (-3.5,-1) 
 -- (0.9,-1) -- (0.9,1) -- 
 (-3.5,1) --
 cycle;
}
}
}
{
\filldraw[draw=black!30!white,very thin,fill=black!40,fill opacity=0.3]%
(140:9mm)
.. controls  +(230:2mm) and (-1,0) ..
(c.east)
.. controls (-1,0) and +(130:2mm) ..
(220:9mm)
arc (220:140:9mm);
}
{
\draw[black!70!white, very thin]
(-3.3,0) -- (-1.9,0);
}
{
\filldraw%
[%
draw=white,
very thin,%
fill=white,%
]
(200:35mm) arc(200:160:35mm)
 -- (-4,1.5) 
 -- (-6,1.5) -- (-6,-1.5) -- 
 (-4,-1.5)
 -- cycle;
\filldraw%
[%
draw=white,
very thin,%
fill=white,%
]
(c) +(20:43mm) arc(20:-20:43mm)
 -- (-3.5,-1.5) 
 -- (1.5,-1.5) -- (1.5,1.5) -- 
 (-3.5,1.5) --
 cycle;
\node at (-2.6,-0.3) {};
}

\node[vertex]
(c) at ( -6, 0) {};
\path (c) node[below left] {$e$};
{
\node[vertex]
(X) at ( 0, 0) {};
\path (X) node[below right] {$x$};
}
}
\end{tikzpicture}
\caption{Restriction to the function $H(x,n)$}\label{figHxn}
\end{figure}
\noindent The following three lemmas provide bounds on the $p$-norms of these functions.

\begin{lem}\label{1} \h There exists some constant $C$ such that for all $x\in X$, $k\leq\frac{n}{4}$ and $n\in\N\setminus\set{0}$,
$\norm{F(x,k,n)}_p^p \leq Cn$.
\m
If, in addition, $d(x,e)\geq 2n$, then $n-3\delta \leq \norm{F(x,k,n)}_p^p$
\end{lem}
\pf The first inequality is obvious as $\abs{F_{x,k,n}}\geq n - 3\delta$. For the second we use Lemma \ref{hypstab} and the bounded geometry of $X$,
\[
\norm{F(x,k,n)}_p^p\ =\ \norm{F(x,k,n)}_1 
\  
\leq 
\  
(2n+1)N(3\delta) 
\  
\leq 
 \  
3N(3\delta)n \leq Cn
\]
completing the proof. \qed

\begin{lem}\label{2} \h If $d(x,y)\leq R$ then $\norm{H(x,n)-H(y,n)}_p\leq 2C(R+1)n^{-\frac{p-1}{p}}$. 
\end{lem}
\pf Choose $x,y\in X$ with $d(x,y)\leq R$. Then as $F_{x,k,n}\subseteq F_{y,k+R,n}$, \vspace{2mm}
\[
\begin{array}{rcl}
\displaystyle\sum_{0\leq k\leq\frac{n}{4}}F(x,k,n)
 &
\leq
 &
\displaystyle\sum_{\frac{n}{4}-R\leq k\leq\frac{n}{4}}F(x,k,n)\ +\ \displaystyle\sum_{0\leq k < \frac{n}{4}-R} F(y,k+R,n)
\m
 &
\leq
 &
\displaystyle\sum_{\frac{n}{4}-R\leq k\leq\frac{n}{4}}F(x,k,n)\ +\ \displaystyle\sum_{0\leq k\leq \frac{n}{4}} F(y,k,n).
\end{array}
\]
Switching $x$ and $y$ in the above argument we conclude that \vspace{2mm}
\[
\begin{array}{rcl}
\norm{H(x,n)-H(y,n)}^p_p
 & 
\leq
 & 
\displaystyle  \frac{\displaystyle 1}{\displaystyle n^p}\displaystyle\sum_{\frac{n}{4}-R\leq k \leq \frac{n}{4}}\norm{F(x,k,n)}^p_p+\norm{F(y,k,n)}^p_p
\m
 &
\leq
 &
\displaystyle \frac{\displaystyle 1}{\displaystyle n^p}(2C(R+1))n \leq (2C(R+1))n^{-(p-1)}.
\end{array}
\]
Notice we have made no assumption that $H(x,n)$, $H(y,n) \neq 0$.\qed

\begin{lem}\label{3} \h $\norm{H(x,n)}^p_p \preceq n$, and whenever $d(x,e)\geq 2n$, $\norm{H(x,n)}^p_p \asymp n$.
\end{lem}
\pf For the lower bound on $\norm{H(x,n)}_p^p$ we notice that given any fixed geodesic $\geo{g}\in\setg{x}{e}$, $\geo{g}([n,2n])\subseteq F(x,k,n)$ for all $k$.
\m
The function $H(x,n)$ takes value at least $\frac{1}{4}$ on at least $n-3\delta$ points, so the lower bound is justified.
\m
As an upper bound, 
\[
\norm{H(x,n)}_p\ \leq\ n^{-1} \displaystyle\sum_{k\leq\frac{n}{4}}\norm{F(x,k,n)}_p\ \leq\ n^{-1}\left(\frac{n}{4}+1\right)(Cn)^{\frac{1}{p}}\ \preceq\ 
n^{\frac{1}{p}}. 
\]
\vspace{-7mm}

\qed
\m
With these three lemmas we are now in a position to define our embedding $\phi:X\to\bigoplus_n\ell^p(X)$.
\[
 \phi(x)\coloneqq\sum_{n\geq 1} \frac{\displaystyle f(2^n)}{\displaystyle 2^{\frac{n}{p}}}H(x,2^n).
\]
To show $\phi$ is Lipschitz, consider $x,y\in X$ with $d(x,y)\leq R$. Then, using Lemma $\ref{2}$:
\[
\begin{array}{rcl} 
\norm{\phi(x)-\phi(y)}^p_p
 &
\preceq
 &
\displaystyle\sum_{n=1}^\infty \frac{f(2^n)^p}{\displaystyle 2^n} \norm{H(x,2^n)-H(y,2^n)}^p_p
\m
 &
\preceq
 &
\displaystyle\sum_{n=1}^\infty \left(\frac{\displaystyle f(2^n)}{\displaystyle 2^n}\right)^p.
\end{array}
\]
As $f$ is concave and has property ($C_p$), $\displaystyle f(2^n)^p \leq 2^{p+1}\sum_{i=2^n+1}^{2^{n+1}}\frac{1}{i}f(i)^p$. Thus
\[
\norm{\phi(x)-\phi(y)}^p_p
 \ 
 \preceq 
 \ 
\displaystyle\sum_{n=1}^\infty \left(\frac{\displaystyle f(2^n)}{\displaystyle 2^n}\right)^p
\ 
\preceq
\ 
\displaystyle\sum_{n=1}^\infty 2^{-np} \sum_{i=2^n+1}^{2^{n+1}}\frac{1}{i}f(i)^p 
 \ 
 \preceq 
 \ 
\displaystyle\sum_{i=1}^\infty \frac{1}{i}\left(\frac{\displaystyle f(i)}{\displaystyle i}\right)^p \preceq 1.
\]
For the lower bound on $\phi$, consider two points $x,y\in X$, with $d(x,y)>12\delta$. We assume, without loss of generality, that $d(x,e)\geq d(y,e)$.
\\
We wish to find a value $k_x$ such that $2^{k_x}\asymp d(x,y)$ and for all $n\in\set{1,2,\dots,k_x}$, the functions $H(x,2^n)$ and $H(y,2^n)$ have disjoint support. Lemma \ref{hypstab} implies that setting
\[
 k_x\coloneqq \lfloor\log_2 \left((x.y)_e - 5\delta\right)\rfloor
\]
suffices, where $(x.y)_e=\frac{1}{2}\left(d(x,e)+d(x,y)-d(y,e)\right)$ is the Gromov product.
\m
Then, by Lemma $\ref{3}$,
\[
\norm{\phi(x)-\phi(y)}_p^p 
 \ 
\succeq 
 \ 
\displaystyle\sum_{n=1}^{k_x} \frac{\displaystyle f(2^{n})^p}{\displaystyle 2^n} \norm{H(x,2^n)}_p^p
 \ 
\succeq
 \ 
\displaystyle\sum_{n=1}^{k_x} f(2^{n})^p
 \ 
\succeq f(d(x,y))^p.
\]
The final step is due to the fact that $f$ is concave and has property ($C_p$).\qed
\m
One may wish to compare the conclusions of Lemmas \ref{1} and \ref{2} with the definition of quantitative property (A) given in \cite{Te08} as another possible method of deducing such an embedding satisfies the conclusion of Theorem \ref{hyp'}.
\m
The next section is independent of the current one, and deals with embeddings of tree-graded spaces. The two approaches are then combined in the final section.

\section{Tree-graded spaces} \label{Freepf}
\noindent During this section, we prove that the compression exponent of amalgamated products and HNN extensions over finite groups depend only on the compression of the initial groups. Specifically we prove that any tree-graded graph can be metrically `decomposed' into a collection of pieces and an underlying tree. Embeddings of tree-graded spaces are found by embedding the two `components' separately, in such a way as to preserve the metric of the original space.
\m
We begin with the definition of a tree-graded graph.

\begin{defn}\label{tgdef} \h Let $\Gamma=(V(\Gamma),E(\Gamma))$ be a connected simplicial graph. We say $\Gamma$ is \emph{tree-graded} \tu{(}in the sense of Dru\c{t}u and Sapir, \tu{\cite{DS05})}, with respect to a collection of non-empty connected subgraphs $\mathcal{P}:=\set{\Gamma_i}_{i\in I}$ if the following properties the satisfied.
\begin{enumerate}
\item Every vertex and every simple loop of $\Gamma$ is contained in some $\Gamma_i$.
\item If $i\neq j$, then $\Gamma_i\not\subseteq \Gamma_j$ and $\abs{V(\Gamma_i)\cap V(\Gamma_j)}\leq 1$.
\end{enumerate}
\end{defn}
\noindent In particular, given two finitely generated groups $A$ and $B$ with finite generating sets $S_A$ and $S_B$ respectively, the Cayley graph Cay($A*B$,$S_A\sqcup S_B$) is tree-graded with respect to the collection of pieces given by the cosets of $A$ and $B$ in $G$. 
\m
The following figure illustrates this definition.
\begin{figure}[H]
\centering
\begin{tikzpicture}[xscale=0.7, yscale=0.35, vertex/.style={draw,fill,circle,inner sep=0.3mm}]
\node[vertex]
(c) at ( -1, 0) {};
{
\shade[xslant=0.5,bottom color=gray!70,top color=gray!10]
    (-6.1,-3.1) rectangle +(12.2,6.2);
\draw[xslant=0.5,color=gray!80, very thin] (-6.1,-3.1) grid +(12.2,6.2);
}
{
\shade[xshift=2cm, yshift=-2cm, bottom color=gray!70,top color=gray!10, opacity=0.5]
    (0,0) -- (-3.1,6.1) -- (3.1,6.1) -- (0,0);
\foreach \r in {-3,-2,-1,0,1,2,3}
{
\draw[xshift=2cm, yshift=-2cm, color=gray!80, very thin, opacity=0.5] (\r,6.1) -- (0,0);
}

\foreach \s in {1,2,3,4,5,6}
{
\begin{scope}
\clip (2,-2) -- (-1.1,4.1) -- (5.1,4.1) -- (2,-2);
\draw[xshift=2cm, yshift=-2cm, color=gray!80, very thin, opacity=0.5]
(0,0) circle (\s cm);
\end{scope}
}
}
{
\shade[xshift=2cm, yshift=2cm, bottom color=gray!50,top color=gray!5, opacity=0.5]
    (0,0) -- (-3.1,6.1) -- (3.1,6.1) -- (0,0);
\foreach \r in {-3,-2,-1,0,1,2,3}
{
\draw[xshift=2cm, yshift=2cm, color=gray!60, very thin, opacity=0.5] (\r,6.1) -- (0,0);
}

\foreach \s in {1,2,3,4,5,6}
{
\begin{scope}
\clip (2,2) -- (-1.1,8.1) -- (5.1,8.1) -- (2,2);
\draw[xshift=2cm, yshift=2cm, color=gray!60, very thin, opacity=0.5]
(0,0) circle (\s cm);
\end{scope}
}
}
{
\shade[xshift=-2cm, yshift=2cm, bottom color=gray!50,top color=gray!5, opacity=0.5]
    (0,0) -- (-3.1,6.1) -- (3.1,6.1) -- (0,0);
\foreach \r in {-3,-2,-1,0,1,2,3}
{
\draw[xshift=-2cm, yshift=2cm, color=gray!60, very thin, opacity=0.5] (\r,6.1) -- (0,0);
}

\foreach \s in {1,2,3,4,5,6}
{
\begin{scope}
\clip (-2,2) -- (-5.1,8.1) -- (1.1,8.1) -- (-2,2);
\draw[xshift=-2cm, yshift=2cm, color=gray!60, very thin, opacity=0.5]
(0,0) circle (\s cm);
\end{scope}
}
}

{
\shade[bottom color=gray!60,top color=gray!8, opacity=0.5]
    (0,0) -- (-3.1,6.1) -- (3.1,6.1) -- (0,0);
\foreach \r in {-3,-2,-1,0,1,2,3}
{
\draw[color=gray!70, very thin, opacity=0.5] (\r,6.1) -- (0,0);
}

\foreach \s in {1,2,3,4,5,6}
{
\begin{scope}
\clip (0,0) -- (-3.1,6.1) -- (3.1,6.1) -- (0,0);
\draw[color=gray!70, very thin, opacity=0.5]
(0,0) circle (\s cm);
\end{scope}
}
}

{
\shade[xshift=2cm, yshift=-2cm, bottom color=gray!70,top color=gray!10, opacity=0.5]
    (0,0) -- (-3.1,6.1) -- (3.1,6.1) -- (0,0);
\foreach \r in {-3,-2,-1,0,1,2,3}
{
\draw[xshift=2cm, yshift=-2cm, color=gray!80, very thin, opacity=0.5] (\r,6.1) -- (0,0);
}

\foreach \s in {1,2,3,4,5,6}
{
\begin{scope}
\clip (2,-2) -- (-1.1,4.1) -- (5.1,4.1) -- (2,-2);
\draw[xshift=2cm, yshift=-2cm, color=gray!80, very thin, opacity=0.5]
(0,0) circle (\s cm);
\end{scope}
}
}

{
\shade[xshift=0cm, yshift=5cm, bottom color=gray!40,top color=white, opacity=0.5]
    (0,0) -- (-3.1,6.1) -- (3.1,6.1) -- (0,0);
\foreach \r in {-3,-2,-1,0,1,2,3}
{
\draw[xshift=0cm, yshift=5cm, color=gray!80, very thin, opacity=0.5] (\r,6.1) -- (0,0);
}

\foreach \s in {1,2,3,4,5,6}
{
\begin{scope}[xshift=0cm, yshift=5cm]
\clip (0,0) -- (-3.1,6.1) -- (3.1,6.1) -- (0,0);
\draw[color=gray!80, very thin, opacity=0.5]
(0,0) circle (\s cm);
\end{scope}
}
}

{
\shade[xshift=-4cm, yshift=0cm, bottom color=gray!50,top color=gray!5, opacity=0.5]
    (0,0) -- (-3.1,6.1) -- (3.1,6.1) -- (0,0);
\foreach \r in {-3,-2,-1,0,1,2,3}
{
\draw[xshift=-4cm, yshift=0cm, color=gray!60, very thin, opacity=0.5] (\r,6.1) -- (0,0);
}

\foreach \s in {1,2,3,4,5,6}
{
\begin{scope}[xshift=-4cm]
\clip (0,0) -- (-3.1,6.1) -- (3.1,6.1) -- (0,0);
\draw[color=gray!60, very thin, opacity=0.5]
(0,0) circle (\s cm);
\end{scope}
}
}

{
\shade[xshift=-6cm, yshift=-2cm, bottom color=gray!50,top color=gray!5, opacity=0.5]
    (0,0) -- (-3.1,6.1) -- (3.1,6.1) -- (0,0);
\foreach \r in {-3,-2,-1,0,1,2,3}
{
\draw[xshift=-6cm, yshift=-2cm, color=gray!60, very thin, opacity=0.5] (\r,6.1) -- (0,0);
}

\foreach \s in {1,2,3,4,5,6}
{
\begin{scope}[xshift=-6cm, yshift=-2cm]
\clip (0,0) -- (-3.1,6.1) -- (3.1,6.1) -- (0,0);
\draw[color=gray!60, very thin, opacity=0.5]
(0,0) circle (\s cm);
\end{scope}
}
}

{
\shade[xshift=-6cm, yshift=1cm, bottom color=gray!50,top color=gray!5, opacity=0.5]
    (0,0) -- (-3.1,6.1) -- (3.1,6.1) -- (0,0);
\foreach \r in {-3,-2,-1,0,1,2,3}
{
\draw[xshift=-6cm, yshift=1cm, color=gray!60, very thin, opacity=0.5] (\r,6.1) -- (0,0);
}

\foreach \s in {1,2,3,4,5,6}
{
\begin{scope}[xshift=-6cm, yshift=1cm]
\clip (0,0) -- (-3.1,6.1) -- (3.1,6.1) -- (0,0);
\draw[color=gray!60, very thin, opacity=0.5]
(0,0) circle (\s cm);
\end{scope}
}
}

{
\shade[xshift=5cm, yshift=2cm, bottom color=gray!50,top color=gray!5, opacity=0.5]
    (0,0) -- (-3.1,6.1) -- (3.1,6.1) -- (0,0);
\foreach \r in {-3,-2,-1,0,1,2,3}
{
\draw[xshift=5cm, yshift=2cm, color=gray!60, very thin, opacity=0.5] (\r,6.1) -- (0,0);
}

\foreach \s in {1,2,3,4,5,6}
{
\begin{scope}[xshift=5cm, yshift=2cm]
\clip (0,0) -- (-3.1,6.1) -- (3.1,6.1) -- (0,0);
\draw[color=gray!60, very thin, opacity=0.5]
(0,0) circle (\s cm);
\end{scope}
}
}

\end{tikzpicture}

\caption{A tree-graded space}\label{figtgdefn}
\end{figure}
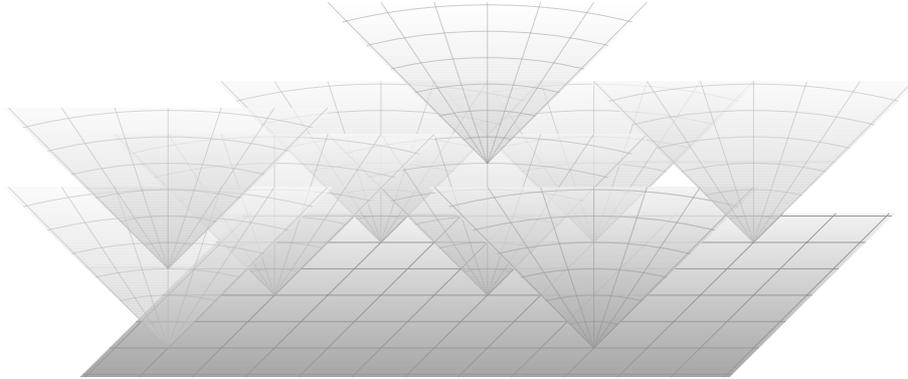
\noindent The main result of this section is

\begin{thm} \label{TG'} \h Let $\Gamma$ be a simplicial graph which is tree-graded with respect to the collection of pieces $\mathcal{P}=\{\Gamma_i\}_{i\in I}$. Suppose we are given a concave function $\rho':\N\to\R_{\geq 0}$ and a collection of coarse embeddings of pieces $\psi_i:\Gamma_i\to\ell^p(X_i)$ such that for all $x,y\in\Gamma_i$,
\[
 \rho'(d_\Gamma(x,y))\leq\norm{\psi_i(x)-\psi_i(y)}_p\leq d_\Gamma(x,y).
\]
If $p=1$, then there is a coarse embedding $\phi$ of $\Gamma$ into an $\ell^1$ space with
\[
 \rho'(d_\Gamma(x,y))\preceq\norm{\phi_i(x)-\phi_i(y)}_p\preceq d_\Gamma(x,y).
\]
For $p>1$, given any function $f:\N\to\R_{\geq 0}$ satisfying \tu{($C_p$)} there is a coarse embedding $\phi$ of $\Gamma$ into an $\ell^p$ space with
\[
 \rho(d_\Gamma(x,y))\preceq\norm{\phi_i(x)-\phi_i(y)}_p\preceq d_\Gamma(x,y),
\]
where $\rho(n)=\min\set{\rho'(n),f(n)}$.
\end{thm}
\noindent We will provide a detailed set-up of the problem which makes the final proof rather short.
\m
Notice that we have made no assumption that $\Gamma$ has any finiteness property, in particular, we have not assumed $I$ is countable. This means that to be absolutely technically accurate we are viewing $\ell^p(X)$ (for this section only) as the space of all countably supported functions $f:X\to\R$ with $\norm{f}^p_p\coloneqq\sum_{i\in I} \abs{f(i)}^p<\infty$. For our purposes, finitely supported functions will suffice.
\m
Theorem \ref{TG'} allows us to calculate compression exponents of certain groups. The list of immediate consequences below is exhaustive by Stallings' Theorem \cite{St68,St71}.

\begin{cor}\label{APandHNN} \h Let $G$ and $H$ be finitely generated groups and let $F$ be a finite subgroup of $G$ and $H$. For all $p\geq 1$,
\begin{enumerate}
 \item $\alpha_p^*(G*_FH)=\min\set{\alpha_p^*(G),\ \alpha_p^*(H)}$ and
 \item $\alpha_p^*(\textup{HNN}(G,F,\theta))=\alpha_p^*(G)$.
\end{enumerate}
\end{cor}
\pf It is an obvious consequence of Theorem \ref{TG'} that (i) holds whenever $F$ is the trivial group.
\m
If $G$ and $H$ are both finite, then $G*_FH$ and $\tu{HNN}(G,F)$ are both hyperbolic by \cite{Gr87}, so the result holds. If $G$ is infinite then $\tu{HNN}(G,F)$ is quasi-isometric to $G*\mathbb{Z}$, similarly, if at least one of $G$ or $H$ is infinite then $G*_FH$ is quasi-isometric to $G*H$, \cite{PW02}. $\Z$ quasi-isometrically embeds into any non-trivial $\ell^p$ space, so we may choose $\rho'(n)=n$, completing the result. \qed
\m
The following lemma is a key tool in this section.

\begin{lem}\label{tggeods}\h Let $\Gamma$ be a tree-graded graph with respect to the collection of subgraphs $\mathcal{P}:=\{\Gamma_i\}_{i\in I}$ and fix a basepoint $e\in V(\Gamma)$.
 For each $x\in V(\Gamma)$ there are finite sets $I_x=\set{i_0,i_1,\dots,i_k}$ and $\set{x_0,\dots,x_k=x}$ such that any geodesic $\geo{g}\in\setg{e}{x}$ can be decomposed into a concatenation of subgeodesics
 \[
  \geo{g} = \geo{g_0}\,\overline{\geo{g_0}}\,\geo{g_1}\,\dots \,\overline{\geo{g_{k-1}}}\,\geo{g_k},
 \]
such that for each $j$,
 \begin{enumerate}
  \item $\geo{g_j}\subseteq \Gamma_j$, where $\iota(\geo{g_j})=e_j$ is the unique vertex of $\Gamma_j$ at minimal distance from $e_0=e$ and $\tau(\geo{g_j})=x_j$,
  \item $\abs{\overline{\geo{g_j}}}\leq 1$.
 \end{enumerate}
\end{lem}
\pf $\geo{g}$ intersects a finite number of pieces, $\setcon{\Gamma_i}{i\in I_x}$ and $\geo{g_i}\coloneqq \geo{g}\cap \Gamma_i$ is connected, so is a sub-geodesic. $\geo{g_i}$ and $\geo{g_{i+1}}$ have at most one vertex in common, by Definition \ref{tgdef}(ii). If they are disjoint, then there is a path of length $1$ between them (otherwise $\Gamma$ must meet another piece, or $I_x$ has been ordered incorrectly). This then gives a decomposition satisfying (i) and (ii).
\m
Now take the above decomposition of any other geodesic in $\setg{e}{x}$. Any counterexample to this lemma forces there to be a simple loop which is not contained in a single piece. \qed
\m
In what follows we will assume, by applying a translation if necessary, that $\psi_i(e_i)=0$ for all $i\in I$.
\m 
Returning to the free product example $G=A*B$, if we are given some word $x=a_0b_1a_2\dots a_{n-1}b_n$, where $a_i\in A\setminus\set{1_A}$ and $b_i\in B\setminus\set{1_B}$, each $\overline{\geo{g_j}}$ is just the vertex $x_j=a_0b_1\dots c_j$, where $c\in\set{a,b}$ depends only on the parity of $j$. Let us assume $c=a$, the other case is similar. Each $\geo{g_j}$ is a geodesic in the coset $x_{j-1}A$ from $x_{j-1}$ to $x_j$, although there may be many choices of such geodesic, they have the same end points.

\begin{defn}\label{defnedisttree}\h Let $\Gamma$ be a tree-graded graph with pieces $\setcon{\Gamma_i}{i\in I}$ and let $e\in V(\Gamma)$. The \emph{$e$-distance tree} of $\Gamma$ is the simplicial tree $T^e_\Gamma$, obtained from $\Gamma$ by identifying vertices under the relation $x\sim y$ if and only if $x,y\in\Gamma_i$ for some $i\in I$ and $d(e_i,x)=d(e_i,y)$.
\end{defn}
\noindent This is not analogous to the coned-off graph of a relatively hyperbolic group, where we imagine collapsing each $\Gamma_i$ to a point, instead we are projecting them to rays.
\m
This construction yields a non-trivial tree even in simple situations like the free product of two groups, where we obtain a tree which does not have bounded geometry, though it is locally finite. What matters here is that geodesics with $e$ as one end vertex inject isometrically into $T^e_\Gamma$ under the obvious graph quotient. This motivates the name $e$-distance tree, as that is precisely what it preserves. We will drop the $e$ from the notation as we have already prescribed a fixed basepoint.
\m
We now introduce a new metric $d'$ on $V(\Gamma)$ which splits distances into a ``tree part'' and a ``pieces part''. $d'$ is the sum of two pseudo-metrics $\sigma_T$ and $\sigma_I$ where $\sigma_T(x,y)$ is the distance between the projections of $x$ and $y$ on the tree $T_\Gamma$ and $\sigma_I(x,y)=\sum_i d(x'_i,y'_i)$ where $x'_i=x_i$ if $i\in I_x$ and $e_i$ otherwise, so $d(x'_i,y'_i)=0$ for all but finitely many values of $i$. 
\m
Checking that $d'$ is a metric only requires showing that $d'(x,y)=0$ implies $x=y$. Assume $d'(x,y)=0$, which ensures $\phi_T(x,y)=0$ and $d(x'_i,y'_i)=0$ for all $i$. The first of these implies that $x$ and $y$ lie in a common piece $\Gamma_j$, and therefore, by definition, $x'_j=x$ and $y'_j=y$. Thus $d(x,y)=d(x'_j,y'_j)=0$ and, as $d$ is a metric, $x=y$.
\m
The following lemma greatly reduces the workload of proving Theorem \ref{TG'}. 

\begin{lem}\label{newmetric}\h Let $d$ be the shortest path metric on $\Gamma$. $d$ and $d'$ are $2$ bi-Lipschitz.
\end{lem}
\pf It is clear that $\sigma_I(x,y)$ and $\sigma_T(x,y)$ are bounded from above by $d(x,y)$, so $d'(x,y)\leq 2d(x,y)$.
\m
For the other bound suppose that for all $i$, $d(x'_i,y'_i)< \frac{1}{2}d(x,y)$, if this is not the case then the result is clear. We show $\sigma_T(x,y)\geq \frac{1}{2}d(x,y)$.
\m
Fix $j\in I_x\cap I_y=\set{i_0,\dots,j=i_l}$ with $d(e,e_j)$ maximal. As $0<d(x'_j,y'_j)< \frac{1}{2}d(x,y)$, in the tree $T_\Gamma$ the images of geodesics in $\setg{x_j}{x}$ and $\setg{y_j,y}$ have at most one vertex in common. 
\\
Moreover, every geodesic in $\setg{e}{x}$ (resp. $\setg{e}{y}$) meets $x_j$ (resp. $y_j$) by Lemma \ref{tggeods}, so $\sigma_T(x,y)\geq d(x,x_j)+d(y,y_j) \geq \frac{1}{2}d(x,y)$. \qed
\m
{\bf Proof of Theorem \ref{TG'}:} \h Results of Tessera and Gal \cite{Te11,Gal} imply that for every $p>1$ and every function $f$ with property ($C_p$) there is a map $\phi_T:V(\Gamma)\to\ell^p(V(T_\Gamma))$ such that
\begin{eqnarray}
 f(\sigma_T(x,y)) \preceq \norm{\phi_T(x)-\phi_T(y)}_p \preceq \sigma_T(x,y).\label{tgphiT}
\end{eqnarray}
While for $p=1$, the above works with $f(n)=n$.
\m
By Lemma \ref{newmetric}, this embedding satisfies the conclusion of Theorem \ref{TG'} whenever $d(x'_i,y'_i)< \frac{1}{2}d(x,y)$ for all $i$. To deal with the other situation, we consider the map $\phi_I: V(\Gamma) \to \bigoplus_{i\in I} \ell^p(V(\Gamma_i))$ given by
\[
 \phi_I(x) = \sum_{i\in I} \psi_i(x'_i).
\]
By definition, $\norm{\phi_I(x)-\phi_I(y)}_p \leq \sigma_I(x,y)$, so this is also Lipschitz. Moreover, given $x,y\in V(\Gamma)$ such that there exists some $i$ with $d(x'_i,y'_i)< \frac{1}{2}d(x,y)$, then
\[
 \norm{\phi_I(x)-\phi_I(y)}_p \succeq \rho(d(x'_i,y'_i)) \succeq \rho(d(x,y)).
\]
Hence the theorem follows for the embedding $\phi=\phi_I + \phi_T$. \qed
\m
However, the embedding $\phi_T$ does not have a good analogue in the relatively hyperbolic case, and as the role of this section is to motivate the methods employed there, we now present another map $\phi_T:V(\Gamma)\to\ell^p(\setcon{e_i}{i\in I})$.
\begin{eqnarray}
(\phi_T(x))(e_k) = 
\left\{
\begin{array}{cl} 
f(d(e_k,x))\left(\displaystyle\frac{d(e_k,e_{k+1})}{d(e_k,x)} \right)^\frac{1}{p}
&
\tu{if}\ k\in I_x\ \tu{and}\ e_k\neq x,\label{phiTweights}
\m
0
&
\tu{otherwise.}
\end{array}
\right.
\end{eqnarray}
\noindent This only works for functions $f$ satisfying property ($C^c_p$), and hence only when $p>1$, as it requires Lemma \ref{ubound} to prove that it is Lipschitz and Lemma \ref{lbound} to provide a lower bound as in Equation (\ref{tgphiT}). We will not give full details here as they are presented in greater generality later.
\m
We now adapt the methods presented so far to the relatively hyperbolic situation.

\section{Relatively hyperbolic spaces}\label{mainthm}
\noindent In this section we mimic the ideas of the previous section, splitting our embedding into two pieces $\phi^s$ and $\phi^l$ which perform the same functions as $\phi_T$ and $\phi_I$ in Section \ref{Freepf} respectively. We also require the averaging techniques from Section \ref{hyp1}. In most of what follows, trying to directly embed an amalgamated product $A*_CB$ where $C$ is finite and non-trivial is sufficient to see why the proof of Theorem \ref{atg'} is so much more intricate than anything previously undertaken in the paper.
\m
We begin with an outline of various properties which are necessary for the method of proof we give. This is a highly notation-heavy process, something we attempt to offset using various informative figures. Once this is complete we formally describe (Definition \ref{keydef}) the collection of simplicial graphs we are considering and prove that asymptotically tree-graded simplicial graphs with bounded geometry satisfy this definition (Proposition \ref{rhinX}).
\m
We assume in what follows that $X$ is the $0$-skeleton of a simplicial graph with bounded geometry equipped with the shortest path metric $d$, $e\in X$ is a fixed basepoint and
\[
 \mathcal{A}=\setcon{A_i}{i\in I}
\]
is a countable collection of subsets - which we will refer to as \emph{pieces} - of $X$ which cover, i.e. $\displaystyle\bigcup_{i\in I} A_i = X$.
\m
In the case of a relatively hyperbolic group, we add the trivial subgroup to the list of peripheral subgroups and take $\mathcal{A}$ to be the set of $K$-neighbourhoods of cosets of these peripheral subgroups, for some fixed $K\geq 1$.
\m
Firstly, we impose some restrictions on geodesics in $X$, in what follows we assume that $\geo{g}\in\setg{x}{e}$, where $e$ is the fixed basepoint of $X$ and $x\in X\setminus\set{e}$:
\m
{\bf Interactions between geodesics and pieces:} 
\m 
We define the $i$-\emph{domain} of $\geo{g}$, denoted $\geo{g}|_i$, to be the convex hull (in $\geo{g}$) of $\geo{g}\cap A_i$ and considering $\geo{g}|_i$ as a directed path, we denote its initial vertex $\iota(\geo{g}|_i)$ by $\geo{g}|_i^+$ and its terminal vertex $\tau(\geo{g}|_i)$ by $\geo{g}|_i^-$.
\begin{figure}[H]
\centering
\begin{tikzpicture}[yscale=0.9,xscale=0.9, vertex/.style={draw,fill,circle,inner sep=0.3mm}]
\clip (-6.5,-0.65) rectangle (6.5,0.8);

\node[vertex]
(E) at (-6, 0) {};
\path (E) node[left] {$e$};

\node[vertex]
(X) at (6, 0) {};
\path (X) node[right] {$x$};

\draw[black!50!white, very thin] (X) -- (E);

\fill[black!30!white, opacity = 0.4]
           (0,0) circle (3cm)
           (3.5,0) circle (0.3cm);
           
\path (0.3,0.1) node[above] {$A_i$};

\node[vertex]
(-) at (-3, 0) {};
\path (-) node[below] {$\geo{g}|_i^-$};

\node[vertex]
(+) at (3.8, 0) {};
\path (+) node[below] {$\geo{g}|_i^+$};[The proof of Theorem $\ref{atg}$]

\draw[black, thick] (-) -- (+);

\path (0,0) node[below] {$\geo{g}|_i$};

\end{tikzpicture}

\caption{The $i$-domain of a geodesic}\label{figidomain}
\end{figure}
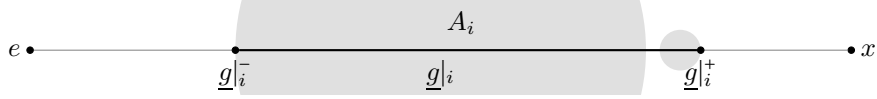
\noindent The $i$-\emph{length} of $\geo{g}$ is defined to be $l_i(\geo{g})\coloneqq d\left(\geo{g}|_i^-,\geo{g}|_i^+\right)+1$.
\m
We require a result analogous to Lemma \ref{hypstab} to ensure an averaging technique like the one in Section \ref{hyp1} can be applied. Hence we impose the following two conditions on $X$ and $\mathcal{A}$.
\begin{enumerate}[label={(C\arabic*)}]
 \item There exists a constant $K>0$ such that for every $i\in I$ and any two geodesics $\geo{g_1}$ and $\geo{g_2}$ which intersect $A_i$, 
 \[
  d\left(\geo{g_1}|_i^-,\geo{g_2}|_i^-\right)\leq K.
 \]
 \item There exists a constant $K>0$ such that for every $i\in I$, every pair of points $x,y$ with $d(x,y)\leq\max\set{\frac{d(x,A_i)}{4},1}$ and any two geodesics $\geo{g_x}\in\setg{x}{e}$ and $\geo{g_y}\in\setg{y}{e}$, which intersect $A_i$, 
 \[
 d(\geo{g_x}|_i^+,\geo{g_y}|_i^+)\leq K.
 \] 
 Moreover, if $\geo{g_y}\cap A_i = \emptyset$, then $l_i(\geo{g_x})\leq K$.
\end{enumerate}
It is a simple consequence of these two conditions that $\abs{l_i(\geo{g_x})-l_i(\geo{g_y})}\leq 2K$ under the hypotheses of condition (C2).
\m
We illustrate the above conditions with the following figure.
\begin{figure}[H]
 \centering
  \begin{tikzpicture}[yscale=2.55,xscale=1.7, 
vertex/.style={draw,fill,circle,inner sep=0.3mm}]
\clip (-6.5,-0.5) rectangle (0.5,0.5);

\fill[black!30!white, opacity = 0.4]
           (-3.3,0) circle (1.3cm);

\node[vertex]
(c) at ( -6, 0) {};

\filldraw[draw=black!50!white, fill=black!35!white, fill opacity=0.7] (0,0) circle (5mm);

\begin{scope}
\clip (0, 0.5) -- (0,-0.5) -- (-3.3,-1.3) arc(-90:90:1.3cm) -- (0,0.5);

\fill%
[%
black!25!white,%
opacity=0.6
]
(120:5 mm)
.. controls +(210:6.5 mm) and (-1.8+0.5,0) ..
(c.east)
.. controls (-1.8+0.5,0) and +(150:6.5 mm) ..
(240:5 mm)
arc (240:120:5 mm);

\draw[black]
(90:4.5 mm)
.. controls +(220:6.5 mm) and (-1.8+0.5,0) ..
(c.east);
\end{scope}

\begin{scope}
\clip (-3.3,0) circle (1.3cm);

\fill%
[%
black!25!white,%
opacity=0.6
]
(-1.92,0) -- (-3.3,1) -- (-4.65,0) -- (-3.3,-1) -- (-1.92,0);
\end{scope}

\begin{scope}
\clip (-6, 0) -- (-3.3,-1.3) arc(-90:90:1.3cm) -- (-6,0);
\draw[black]
(-1.855,0) -- (-3.3, 0.4) -- (-4.6,0.02) -- (-6,0);
\end{scope}

\node[vertex]
(c) at ( -6, 0) {};
\path (c) node[below left] {$e$};

\node[vertex]
(x) at ( 0, 0) {};
\path (x) node[below] {$x$};

\path (0,0.03) node[above] {$B(x;\frac{n}{4})$};
\path (-3.3,0.37) node[below] {$\geo{g_y}$};

\draw[black] (c) -- (x);
\path (-0.65,0) node[below] {$\geo{g_x}$};

\node[vertex]
(y) at ( 90:4.5mm) {};
\path (y) node[below right] {$y$};

\draw[black, <->, thin] (-2,-0.24) -- (0,-0.24);
\path (-1,-0.24) node[below] {$n$};

\path (-3.3,-0.3) node[] {$A_i$};
\end{tikzpicture}
\caption{A relative version of Lemma \ref{hypstab}}\label{figrelhypstab}
\end{figure}
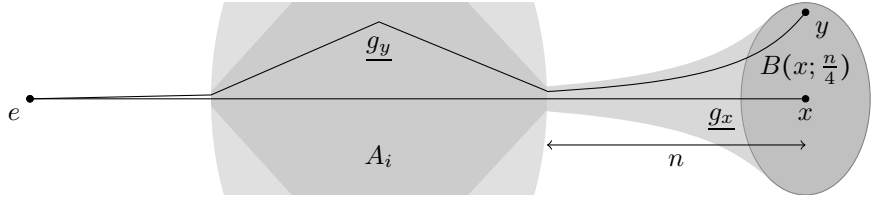
\noindent For relatively hyperbolic groups, both these conditions follow from the bounded coset penetration property, originally defined by Farb, \cite{Fa98}.
\m
A little more notation is required here for the next stages. 
\m
We will also need various finiteness conditions to ensure Lipschitz upper bounds for the averaging procedures.
\m
{\bf Finiteness conditions:}
\m
Figure \ref{figrelhypstab} illustrates a key finiteness issue we will have to address, namely, if we consider the collection of all geodesics from a ball of radius $k$ around $x$ to the origin
\[
\geo{G_{x,k}} \coloneqq \bigcup_{d(x,y)\leq k} \setg{y}{e},
\]
then, within pieces where such geodesics have large $i$-domains, this collection of geodesics may meet a large number of other pieces. As a result, any embedding we make taking the form of those in Sections \ref{hyp1} and \ref{Freepf} will fail to be Lipschitz in general. We remedy this by discounting pieces whose $i$-domains lie `deep inside' some large $j$-domain
\m
For a given $x\in X\setminus\set{e}$ we restrict the collection of considered pieces in the following way:
\m
Given a collection of geodesics $\geo{G}$, and a constant $K\geq 1$ we define the \emph{$i$-boundary of $\geo{G}$}, $\partial^K_i(\geo{G})$ to be the set of $\geo{g}|_i^+$ which satisfy the following condition.
\m
For all $\geo{g'}\in\geo{G}$ and any $j\in I$ with $l_j(\geo{g'})\geq 5K$,
\begin{equation}
 d(e,\geo{g}|_i^+)\notin \left[d(e,\geo{g'}|_j^-)+2K, d(e,\geo{g'}|_j^+)-2K \right].\label{eqIxK}
\end{equation}
This unpleasant-looking restriction discounts those $i$-domains which lie `deep' inside some $j$-domain. In Figure $6$ one can imagine that by taking $A_i$ to be a copy of $\Z^2$ in some free product, the collection of all geodesics from a fixed $x$ to $e$ with $i$-domain of length $n$ meet around $n^2$ different pieces while inside $A_i$. This restriction - equation (\ref{eqIxK}) - throws out all but a uniformly bounded number of such domains and is crucial in showing that the two maps we will construct are Lipschitz, which we do in Lemmas \ref{phi^sLip} and \ref{phi^lLip}.
\m
Then when we consider the analogue of the ``hyperbolic trumpets'' considered in Section \ref{hyp1}, we collect all suitable boundary points of relevant $i$-domains into the set
\[
 \partial^K_i(x)\coloneqq \bigcup_{k\leq d(x,A_i)/4}\partial^K_i\left(\geo{G_{x,k}}\right).
\]
We define the set of \emph{$x$-relevant $i$-domains} - $I_x(K)$ - to be the set of $i\in I$ which are crossed by some geodesic in a ``trumpet'' around $x$, but not too close to the basepoint. More formally, we require that
\[
 \partial^K_i(x)\neq\emptyset \tu{ and } d(e,\partial^K_i(x)) \geq 3K.
\]
It is far simpler - and thus very tempting - to consider only pieces with a sufficiently large $i$-domain. However, simply considering the situation of a hyperbolic (but not free) group as hyperbolic relative to the trivial group, we obtain an empty collection of pieces.
\m
As in the hyperbolic case, we dismiss any piece which is too close to the basepoint. 
\m
The subset $I'_x(K)\subseteq I_x(K)$ consists of those $i\in I_x(K)$ such that $x\notin A_i$.
\m
{\bf Technical point:} \h Unlike the hyperbolic situation (cf. Lemma \ref{2}), it is not in general true that $\partial^K_i(\geo{G_{x,k}})\subseteq \partial^K_i(\geo{G_{y,k+R}})$ whenever $d(x,y)\leq R$. However, we do have the following.

\begin{lem}\label{lemnxi} \h For each $a\in X$ define $n_{x,i}(a)\coloneqq\abs{\setcon{k}{a\in \partial^K_i(\geo{G_{x,k}})}}$. Then, for any two points $x,y\in X$ with $d(x,y)\leq R$,
\[
 \abs{n_{x,i}(a)-n_{y,i}(a)}\leq\abs{\setcon{k}{a\in \partial^K_i(\geo{G_{x,k}})\symd \partial^K_i(\geo{G_{y,k}})}}\leq 4R,
\]
where $\symd$ denotes the symmetric difference $A\symd B\coloneqq (A\setminus B) \sqcup (B\setminus A)$.
\end{lem}
\pf We bound the number of possible values of $k$ such that 
\[
 a\in \partial^K_i(\geo{G_{x,k}})\setminus \partial^K_i(\geo{G_{y,k}})
\]
For this to occur, either $a$ is not an end point of any $i$-domain of a geodesic in $\geo{G_{y,k}}$, or it is, but there is some $\geo{g'}\in\geo{G_{y,k}}$ and some $j\in I$ such that Equation (\ref{eqIxK}) fails to hold.
\m
In the second situation, $a\not\in \partial^K_i(\geo{G_{x,m}})\cup \partial^K_i(\geo{G_{y,m}})$ for any $m\geq k+R$. 
\m
In the first situation we are only interested in the case where $a\not\in \partial^K_i(\geo{G_{y,k+R}})$, which only occurs when the second situation holds for $k+R$. Thus,
\[
 \abs{\setcon{k}{a\in \partial^K_i(\geo{G_{x,k}})\setminus \partial^K_i(\geo{G_{y,k}})}}\leq 2R,
\]
and the bound on $\abs{n_{x,i}(a)-n_{y,i}(a)}$ follows. \qed
\m
This lemma will be sufficient to develop an analogue of Lemma \ref{2}, once we have the following finiteness properties.

\begin{enumerate}[resume,label={(C\arabic*)}]
 \item There exists a constant $K$ such that $1\leq\abs{\setcon{i}{x\in A_i}}\leq K$ for all $x\in X$ and $\tu{diam}(A_i\cap A_j)\leq K$ for all $i,j\in I$ with $i\neq j$.
 \item There exists a constant $K$ such that $\abs{\setcon{i\in I_x(K)}{d(x,A_i)=t}}\leq K$, for each $t\in\N$.
\end{enumerate}
\noindent We are now ready to define the collection of spaces considered.
\begin{defn}\label{keydef} \h {\bf SPQR relatively hyperbolic graphs}
 \m
 Let $X$ be the $0$-skeleton of a simplicial graph with bounded geometry and a fixed basepoint $e\in X$. Let $\mathcal{A}=\setcon{A_i}{i\in I}$ be a collection of subsets of $X$. The triple $(X,\mathcal{A},e)$ is \emph{SPQR relatively hyperbolic} if it satisfies conditions \tu{(C1)-(C4)} for a fixed constant $K$.
\end{defn}
\noindent Intuitively, what this says is that $X$ is SPQR relatively hyperbolic if it appears to be a relatively hyperbolic graph when all we can understand is geodesics to the basepoint $e$, i.e. all those roads that lead to Rome. We expect that any Cayley graph which is SPQR relatively hyperbolic is actually the Cayley graph of a relatively hyperbolic group.
\m
The main examples of SPQR relatively hyperbolic spaces are given by the following proposition.

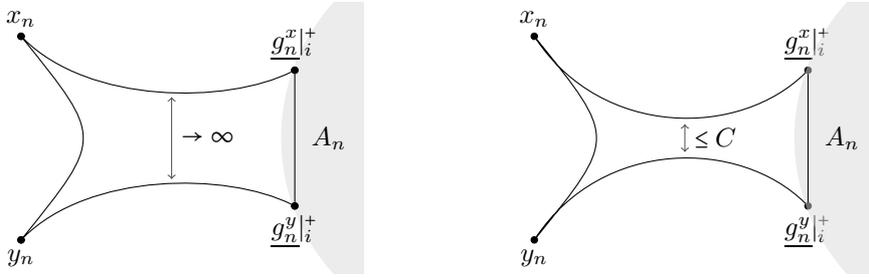
\begin{figure}[H]
 \centering
 \begin{tikzpicture}[yscale=0.9,xscale=0.9, vertex/.style={draw,fill,circle,inner sep=0.3mm}]

\node[vertex]
(X) at (-2, 1.5) {};

\node[vertex]
(T) at (-2, -1.5) {};

\node[vertex]
(Y) at (2, 1) {};

\node[vertex]
(Z) at (2, -1) {};

\draw[black!70!white, very thin, <->]
        (0.2,0.6) -- (0.2,-0.6);
\path (0.2,0) node[right] {$\to\infty$};

\begin{scope}
 \clip (0,-2) rectangle (3,2);
 \fill[black!15!white, opacity = 0.5]
               (4.8,0) circle (3cm);
\end{scope}
\draw[] (X) .. controls (-1,0.5) and (1,0.5) .. (Y)
        (T) .. controls (-1,-0.5) and (1,-0.5) .. (Z)
        (X) .. controls (-0.8,0) .. (T)
        (Y)--(Z);
\path (2.5,0) node[] {$A_n$};

\node[vertex]
(X) at (-2, 1.5) {};
\path (X) node[above] {$x_n$};

\node[vertex]
(T) at (-2, -1.5) {};
\path (T) node[below] {$y_n$};

\node[vertex]
(Y) at (2, 1) {};
\path (Y) node[above] {$\geo{g^x_n}|_i^+$};

\node[vertex]
(Z) at (2, -1) {};
\path (Z) node[below] {$\geo{g^y_n}|_i^+$};

\begin{scope}[xshift=7.5cm]

\node[vertex]
(X) at (-2, 1.5) {};
\path (X) node[above] {$x_n$};

\node[vertex]
(T) at (-2, -1.5) {};
\path (T) node[below] {$y_n$};

\node[vertex]
(Y) at (2, 1) {};
\path (Y) node[above] {$\geo{g^x_n}|_i^+$};

\node[vertex]
(Z) at (2, -1) {};
\path (Z) node[below] {$\geo{g^y_n}|_i^+$};

\draw[black!70!white, very thin, <->]
        (0.2,0.2) -- (0.2,-0.2);
\path (0.2,0) node[right] {$\leq C$};

\begin{scope}
 \clip (0,-2) rectangle (3,2);
 \fill[black!15!white, opacity = 0.5]
               (4.8,0) circle (3cm);
\end{scope}
\path (2.5,0) node[] {$A_n$};

\draw[] (X) .. controls (-1,0) and (1,0) .. (Y)
        (T) .. controls (-1,0) and (1,0) .. (Z)
                (X) .. controls (-0.8,0) .. (T)
                (Y)--(Z);
\end{scope}

\end{tikzpicture}

 \caption{Fat and `pinched' quadrilaterals}\label{figfatquads}
\end{figure}

\begin{prop}\label{rhinX} \h Let $X$ be a simplicial graph of uniformly bounded degree which is asymptotically tree-graded in the sense of \tu{\cite{DS05}}, with set of pieces $\mathcal{A}'$. Then $(X,\mathcal{A},x)$ is SPQR relatively hyperbolic where $\mathcal{A}$ is a set of $M$-neighbourhoods of elements of $\mathcal{A}'$ and $M$-balls around points and $x\in X$ is arbitrary.
\end{prop}
\pf Suppose $X$ is asymptotically tree-graded with respect to a collection of subsets, which we will label $\setcon{A_i}{i\in I}$. Then we set $\mathcal{A}$ to be the collection of all $M$-neighbourhoods of these pieces and all $M$ balls centred at points not lying in some $A_i$, where $M$ is the constant obtained in the proof of the Rips' hyperbolicity of saturations, \cite[$4.27$]{DS05}.
\m
Property (C1) is the conclusion of \cite[$8.14$]{DS05}. If we suppose that (C2) fails, and choose a collection of counterexample triples $(x_n,y_n,A_n)$ such that $d(x_n,y_n)\leq d(x_n,A_n)/4$ and there are geodesics $\geo{g^x_n}\in\setg{x_n}{e}$ and $\geo{g^y_n}\in\setg{y_n}{e}$ such that
 \[ 
  d(\geo{g^x_n}|_i^+,\geo{g^y_n}|_i^+)\geq n.
 \]
Quadrilaterals with (ordered) vertex set $(x_n,\geo{g^x_n}|_i^+,\geo{g^y_n}|_i^+, y_n, x_n)$ are either uniformly pinched - in the sense of Figure \ref{figfatquads} below - or can be capped to produce polygons which are increasingly fat in the sense of \cite[Definition $3.32$]{DS05}.
\noindent If they are increasingly fat, then we obtain a contradiction to \cite[Theorem $4.1$($\alpha_3$)]{DS05} as such polygons are certainly not contained in a uniform neighbourhood of a single piece. However, if they are pinched we obtain a contradiction to \cite[Lemma $4.11$]{DS05}, which can be thought of as a quasi-geodesic version of property (C1). Here we are using the quasi-convexity of pieces in an asymptotically tree-graded space \cite[Lemma $4.3$]{DS05}.
\m
Property (C3) follows verbatim from \cite[$4.1$($\alpha_1$)]{DS05}. Finally, property (C4) follows from the other three properties and the bounded geometry of $X$. Properties (C1) and (C2) ensure that we only need consider points in uniformly bounded neighbourhood of any fixed geodesic $\geo{g_x}\in\setg{x}{e}$ - here the restriction to $i\in I_x(K)$ is crucial. This contains a uniformly bounded number of points at any fixed distance from $x$, each of which lies in only finitely many pieces, by (C3).\qed
\m
In particular, Cayley graphs of relatively hyperbolic groups with finitely generated peripheral subgroups are SPQR relatively hyperbolic, \cite[Appendix A]{DS05}.
\m
Now we present the main theorem of the paper in its most general form.

\begin{thm}\label{atg'} \h \tu{(cf. Theorem \ref{RHgps})}
\m
Let $(X,\mathcal{A},e)$ be SPQR relatively hyperbolic, where $\mathcal{A}=\setcon{A_i}{i\in I}$. Suppose we are provided with maps $\psi_i:A_i\to\ell^p(X_i)$ and a concave function $\rho':\R_{\geq 0}\to\R_{\geq 0}$ such that for all $x,y\in A_i$,
\[
 \rho'(d(x,y))\leq\norm{\psi_i(x)-\psi_i(y)}_p\leq d(x,y).
\]
Then, for all $p>1$ and all functions $f:\N\to\R_{\geq 0}$ with property \tu{($C_p^c$)} there exists a coarse embedding $\phi$ of $X$ into some $\ell^p$ space such that for all $x,y\in X$,
\[
\rho(d(x,y))\preceq\norm{\phi(x)-\phi(y)}_p\preceq d(x,y),
\]
where $\rho(n) = \min\set{\rho'(n),f(n)}$.
\end{thm}
\noindent We define $e_i$ to be some closest point of $A_i$ to $e$, by condition (C3), the diameter of the set of possible choices for $e_i$ is at most $K$. Without loss of generality we may assume $\psi_i(e_i)=0$ for each $i\in I$. As the constant $K$ is now fixed we will write $I_x$ and $I'_x$ in place of $I_x(K)$ and $I'_x(K)$ respectively. Similarly, we drop the $K$ in the notation $\partial_i^K$. 
\m
The proof is now split into three parts, in the first two we introduce two maps from $X$ into $\ell^p$ spaces with the equivalent roles of $\phi_T$ and $\phi_I$ in Section \ref{Freepf} and prove they are Lipschitz. In the third (\ref{pfatg}) we prove their sum satisfies the conclusion of Theorem $\ref{atg'}$.

\subsection{Embedding small pieces}
We construct an embedding with the role of identifying and separating points whose geodesics spend most of their time in pieces they see very little of.
\m
It is important to note that this uses the technique of replacing each piece by a ray like in the construction of the $e$-distance tree (Definition \ref{defnedisttree}). We then use the averaging method from Section \ref{hyp1} to ensure a Lipschitz map. Crucially, in what follows, we can only average a length proportional to the distance between a point $x$ and a piece $A_i$ - in the hyperbolic case we take geodesics of length $n$ at distance approximately $n$ from $x$. For this reason, we define the following capped version of the hyperbolic trumpets $F_{x,k}$ seen in Section \ref{hyp1}.
\m
For each $i\in I'_x$, we define functions $F_i(x,k)\in\ell^p(X)$ as follows:
\[
 F_i(x,k)(y) = \left\{
\begin{array}{ll}
\min\set{d(x,A_i), d_X(y,e_i)+1}^\frac{1}{p}
 &
\tu{if } y\in \partial_i(\geo{G_{x,k}})
\m
0
 &
\tu{otherwise.}
\end{array}\right.
\]
As a necessary shorthand we set $d_{x,i}(y)\coloneqq\min\set{d(x,A_i), d_X(y,e_i)+1}$. We then define 
\[
 H_i(x)=\frac{1}{d(x,A_i)} \sum_{k\leq\frac{d(x,A_i)}{4}} F_i(x,k).
\] 
The following three lemmas (mirroring Lemmas \ref{1}, \ref{2} and \ref{3}) provide useful information on these new objects.

\begin{lem}\label{1'}\h For all $x\in X$, $\geo{g}\in\setg{x}{e}$, $i\in I'_x$ and $k\leq\frac{d(x,A_i)}{4}$,
\[
 d_{x,i}(\geo{g}|_i^+)\leq \norm{F_i(x,k)}_p^p \leq \abs{\partial_i(\geo{G_{x,k}})}\left(d_{x,i}(\geo{g}|_i^+) +K\right) \preceq d_{x,i}(\geo{g}|_i^+).
\]
Moreover, if $\displaystyle\geo{g}|_i^+\in\bigcap_{k\leq d(x,A_i)/4}\partial_i(\geo{G_{x,k}})$, then
\[
 \norm{F_i(x,k)}_p^p \geq d_{x,i}(\geo{g}|_i^+).
\]
\end{lem}
\pf By assumption, $\geo{g}|_i^+\in\partial_i(\geo{G_{x,k}})$ for all $k\leq\frac{d(x,A_i)}{4}$, so the first bound is satisfied.
\m
For the second bound, 
\[
\norm{F_i(x,k)}_p^p \leq \sum_{y\in \partial_i(\geo{G_{x,k}})} d_{x,i}(y) \leq \abs{\partial_i(\geo{G_{x,k}})}(d_{x,i}(\geo{g}|_i^+) + K).
\]
The final inequality holds as $X$ has uniformly bounded valency and the diameter of $\partial_i(\geo{G_{x,k}})$ is at most $K$, by condition (C2).

\begin{lem}\label{2'} \h For all $x\in X$, $i\in I'_x$ and $\geo{g}\in\setg{x}{e}$,
\[
 \norm{H_i(x)}_p^p \leq \abs{\partial_i(x)}\left(d_{x,i}(\geo{g}|_i^+) +K\right). 
\]
Moreover, if $\displaystyle\geo{g}|_i^+\in\bigcap_{k\leq d(x,A_i)/4}\partial_i(\geo{G_{x,k}})$, then
\[
 \norm{H_i(x)}_p^p\geq \frac{1}{4}d_{x,i}(\geo{g}|_i^+).
\]
\end{lem}
\pf The upper bound follows from Lemma \ref{1'}, as $\geo{g}|_i^+\in \partial_i(\geo{G_{x,k}})$ for all $k\leq\frac{d(x,A_i)}{4}$, so
\[
\norm{H_i(x)}_p \leq \frac{\displaystyle 1}{\displaystyle d(x,A_i)} \frac{\displaystyle d(x,A_i)+1}{\displaystyle 4} \norm{F_i\left(x,\frac{d(x,A_i)}{4}\right)}_p.
\] 
For the lower bound, we evaluate the contribution to $\norm{H_i(x)}_p$ coming from the point $\geo{g}|_i^+$:
\[
\norm{H_i(x)}_p \geq \frac{\displaystyle 1}{\displaystyle d(x,A_i)} \frac{\displaystyle d(x,A_i)+1}{\displaystyle 4} d_{x,i}\left(\geo{g}|_i^+\right)^\frac{1}{p}.
\]

\begin{lem}\label{3'} \h There exists some constant $C>0$ such that for all $x,y\in X$ with $d(x,y)\leq 1$, all $\geo{g}\in\setg{x}{e}$ and all $i\in I'_x\cup I'_y$,
\[
\norm{H_i(x)-H_i(y)}^p_p \leq C\frac{d_{x,i}(\geo{g}|_i^+)}{d(x,A_i)^p}.
\]
\end{lem}
\pf We first bound the absolute value of $H_i(x)-H_i(y)$ at some point $a\in \partial_i(x)\cup \partial_i(y)$.
\m
Recall that $n_{z,i}(a):=\abs{\setcon{n\leq \frac{d(z,A_i)}{4}}{a\in \partial_i(\geo{G_{z,k}})}}$. So
\[
\abs{\left(H_i(x)-H_i(y)\right)(a)} = \abs{\frac{\displaystyle n_{x,i}(a)}{\displaystyle d(x,A_i)}d_{x,i}(a)^\frac{1}{p} - \frac{\displaystyle n_{y,i}(a)}{\displaystyle d(y,A_i)}d_{y,i}(a)^\frac{1}{p}}.
\]
By Lemma \ref{lemnxi}, we know that whenever $d(x,y)\leq 1$, $\abs{n_{x,i}(a)-n_{y,i}(a)}\leq 4$, for all $a$.
\m
If $i\in I'_x\setminus I'_y$ then $H_i(y)=0$, and $n_{x,i}(a)\leq 2$. Again, we use the fact that $\abs{\partial_i(x)}$ is uniformly bounded by condition (C2) and the uniformly bounded valency of $X$, so we are done. The case $i\in I'_y\setminus I'_x$ is treated in the same way.
\m
Suppose now that $i\in I'_x\cap I'_y$, so $d(x,A_i),d(y,A_i)\geq 1$. Notice that $d_{x,i}(a)=d_{y,i}(a)$ unless one (or both) are equal to $d(x,A_i)$ (respectively $d(y,A_i)$). Therefore,
\[
\abs{d_{x,i}(a)^\frac{1}{p}-d_{y,i}(a)^\frac{1}{p}} \leq \min\set{d(x,A_i),d(y,A_i)}^{-\frac{p-1}{p}} \leq 2d(x,A_i)^{-\frac{p-1}{p}}.
\]
At this point it is crucial that we capped the lengths considered before defining the embedding.
\m
Finally, combining these observations we have
\[
\abs{\left(H_i(x)-H_i(y)\right)(a)} = \frac{\displaystyle d(y,A_i)n_{x,i}(a)d_{x,i}(a)^\frac{1}{p}-d(x,A_i)n_{y,i}(a)d_{y,i}(a)^\frac{1}{p}}{\displaystyle d(x,A_i)d(y,A_i)}.
\]
By the triangle inequality, we can bound this from above by
\[
 \frac{\displaystyle n_{x,i}d_{x,i}(a)^\frac{1}{p}\abs{d(x,A_i)-d(y,A_i)}}{\displaystyle d(x,A_i)d(y,A_i)} + \frac{\displaystyle d(x,A_i)d_{x,i}(a)^\frac{1}{p}\abs{n_{x,i}(a)-n_{y,i}(a)}}{\displaystyle d(x,A_i)d(y,A_i)}
\]
\[
 + \frac{\displaystyle d(x,A_i)n_{y,i}(a)\abs{d_{y,i}(a)^\frac{1}{p}-d_{x,i}(a)^\frac{1}{p}}}{\displaystyle d(x,A_i)d(y,A_i)}.
\]
Applying all the previous deductions and noticing that $n_{x,i}(a)\leq d(x,A_i)$ we obtain a uniform constant $C'$ such that
\[
\abs{\left(H_i(x)-H_i(y)\right)(a)} \leq C' \frac{\displaystyle d_{x,i}(a)^\frac{1}{p}}{\displaystyle d(x,A_i)}.
\]
Finally, we use condition (C2) again to deduce that $\abs{\partial_i(x)\cup \partial_i(y)}$ is uniformly bounded and the lemma follows. \qed
\m
We are now ready to define the first part of our embedding:
\[
 \phi^{\tu{s}}(x):= \sum_{i\in I'_x} \frac{f(d(x,A_i))}{d(x,A_i)^\frac{1}{p}} H_i(x).
\]
\begin{lem}\label{phi^sLip} \h $\phi^s: X\to \ell^p(X)$ is Lipschitz for all $p> 1$.
\end{lem}
\pf Consider two points $x,y\in X$ with $d(x,y)\leq 1$.
\m
Firstly, suppose $i\in I'_x\setminus I'_y$. Then by Lemma \ref{3'},
\[
 \norm{H_i(x)}_p^p \leq C\frac{d_{x,i}(\geo{g}|_i^+)}{d(x,A_i)^p}, 
\]
for any geodesic $\geo{g}\in\setg{x}{e}$, but by condition (C2), $l_i(\geo{g})$ is uniformly bounded, so
\[
 \norm{H_i(x)}_p^p \preceq \frac{1}{(d(x,A_i)+1)^p}.
\]
The case $i\in I'_y\setminus I'_x$ is treated similarly and as $\abs{d(x,A_i)-d(y,A_i)}\leq 1$,
\[
 \norm{H_i(y)}_p^p \preceq \frac{1}{(d(x,A_i)+1)^p}.
\]
By the triangle inequality, the contribution made to $\norm{\phi^s(x)-\phi^s(y)}$ by those $i\in I'_x\cap I'_y$ is at most
\begin{eqnarray}
\sum_{i\in I'_x\cap I'_y} \frac{f(d(x,A_i))^p}{d(x,A_i)} \norm{H_i(x)-H_i(y)}_p^p \label{phis1}
\end{eqnarray}
\begin{eqnarray}
 + \sum_{i\in I'_x\cap I'_y} \left( \frac{f(d(x,A_i))}{d(x,A_i)^\frac{1}{p}} - \frac{f(d(y,A_i))}{d(y,A_i)^\frac{1}{p}}\right)^p \norm{H_i(x)}_p^p. \label{phis2}
\end{eqnarray}
\m
As $f(n)(n^{\frac{-1}{p}})$ is non-decreasing (cf. Definition \ref{defnccp}) we may use the same argument as in the tree-graded case to deduce that $(\ref{phis2})$ is bounded from above (up to some uniform multiplicative constant) by
\[
 \sum_{i\in I'_x\cap I'_y} \frac{\min\{\tu{diam}(\partial_i(x))+1, d(x,A_i)\}}{d(x,A_i)}\left(\frac{f(d(x,A_i))}{d(x,A_i)} \right)^p.
\]
Also, by Lemma \ref{2'} and the fact that $f$ is concave, $(\ref{phis1})$ is bounded from above (up to some uniform multiplicative constant) by
\[
 \sum_{i\in I'_x\cap I'_y} \frac{\min\{\tu{diam}(\partial_i(x))+1, d(x,A_i)\}}{d(x,A_i)}\left(\frac{f(d(x,A_i))}{d(x,A_i)} \right)^p.
\]
Hence,
\[
 \norm{\phi^{\tu{s}}(x)-\phi^{\tu{s}}(y)}_p^p \preceq \sum_{i\in I'_x\cup I'_y} \frac{\min\set{\tu{diam}(\partial_i(x))+1, d(x,A_i)}}{d(x,A_i)}\left(\frac{f(d(x,A_i))}{d(x,A_i)} \right)^p 
\]
which is uniformly bounded. To see this notice that by conditions (C3) and (C4) we can partition $I'_x\cup I'_y$ into a uniformly bounded number of subsets in such a way that the above sum restricted to any such subset satisfies the hypotheses of Lemma \ref{ubound}. \qed

\subsection{Embedding large pieces}
For the second part of the embedding we make a complementary construction, set the task of identifying long pieces using the existing embeddings of pieces $(\psi_i)_{i\in I}$. The difficulty here is in ensuring the map is Lipschitz. To do this we combine geodesics using the averaging methods from Section \ref{hyp1}, normalised so that each `thick geodesic' has suitable weight.
\m
We define $a_{x,i}\coloneqq\displaystyle \sum_{k=0}^{\frac{d(x,A_i)}{4}} \abs{\partial_i(\geo{G_{x,k}})}$ and define $
k_{x,i}=\min\set{a_{x,i},1+\frac{d(x,A_i)}{4}}$. This will be the normaliser of our thick geodesic.
\m
Recall that we made the convention $\psi_i(e_i)=0$ for each $i\in I$.
\m
We then proceed towards the definition of the second part of the embedding, by defining
\[
  F'_i(x,k)\coloneqq\sum_{a\in \partial_i(\geo{G_{x,k}})} \psi_i(a).
\]
We then normalise using $k_{x,i}$,
\[
  H'_i(x)=\frac{\displaystyle 1}{\displaystyle k_{x,i}}\sum_{k\leq \frac{d(x,A_i)}{4}} F'_i(x,k).
\]
The second part of the embedding $\phi^l:X\to \bigoplus_{i\in I}\ell^p(X_i)$ is defined as
\[
  \phi^l(x)=\sum_{i\in I_x} H'_i(x).
\]
\begin{lem}\label{phi^lLip} \h For all $p>1$, $\phi^l$ is Lipschitz.
\end{lem}
\pf Let $x,y\in X\setminus B(e;3K)$ with $d(x,y)\leq 1$. We show that for each $i\in I_x\cup I_y$, 
\[
  \norm{H'_i(x)-H'_i(y)}_p \leq \frac{C}{d(x,A_i)+1}, 
\]
for some $C>0$ not depending on $i$. This suffices by the finiteness conditions (C3) and (C4).
\m
Initially, suppose $k_{x,i}=a_{x,i}$ and $k_{y,i}=a_{y,i}$. Then notice that the function
\[
 \frac{\displaystyle 1}{\displaystyle k_{x,i}}\sum_{k\leq \frac{d(x,A_i)}{4}} \chi(\partial_i(\geo{G_{x,k}}))
\]
where $\chi(S)$ is the characteristic function of the set $S$, is non-negative and has $\ell^1$ norm exactly $1$, as $k_{x,i}=\displaystyle \sum_{k\leq \frac{d(x,A_i)}{4}} \abs{\partial_i(\geo{G_{x,k}})}$.
\m
Moreover,
\begin{eqnarray}
 \frac{\displaystyle 1}{\displaystyle k_{x,i}}\sum_{k\leq \frac{d(x,A_i)}{4}} \chi(\partial_i(\geo{G_{x,k}})) - \frac{\displaystyle 1}{\displaystyle k_{y,i}}\sum_{k\leq \frac{d(y,A_i)}{4}} \chi(\partial_i(\geo{G_{y,k}})) \label{diffl1}
\end{eqnarray}
has $\ell^1$ norm at most $\frac{C'}{d(x,A_i)}$ for some uniform constant $C'$, and the sum of its entries is $0$.
\m
The second of these claims follows from the fact that this is a difference of non-negative functions of $\ell^1$ norm $1$. Recall that $\abs{n_{x,i}(a)-n_{y_i}(a)}\leq 4$ for all $a$, by Lemma \ref{lemnxi} and the set $\partial_i(\geo{G_{x,k}})\cup \partial_i(\geo{G_{y,k}})$ has uniformly bounded cardinality (independent of $k$). Hence, $\abs{k_{x,i}-k_{y,i}}$ is uniformly bounded by some constant $C''$.
\m
Next, fix any point $a\in X$. The contribution to $(\ref{diffl1})$ coming from $a$ is at most
\[
 \abs{\frac{n_{x,i}(a)}{k_{x,i}} - \frac{n_{y,i}(a)}{k_{y,i}}}.
\]
As $n_{y,i}(a)\leq k_{y,i}$,
\[
\begin{array}{rcl} \displaystyle\abs{\frac{n_{x,i}(a)}{k_{x,i}} - \frac{n_{y,i}(a)}{k_{y,i}}} 
 &
 \leq 
 &
\displaystyle \abs{\frac{n_{x,i}(a)}{k_{x,i}}-\frac{n_{y,i}(a)}{k_{x,i}}} + \abs{\frac{n_{y,i}(a)}{k_{x,i}}-\frac{n_{y,i}(a)}{k_{y,i}}}  
\m
 &
 \leq
 &
\displaystyle \frac{\abs{n_{x,i}(a)-n_{y,i}(a)}}{k_{x,i}} + \frac{n_{y,i}(a)\abs{k_{y_i}-k_{x,i}}}{k_{y,i}k_{x,i}}
\m
 &
 \leq
 &
\displaystyle \frac{1}{k_{x,i}} + \frac{C''}{k_{x,i}} \leq \frac{C'}{d(x,A_i)+1},
\end{array}
\]
with the final step coming from the fact that $k_{x,i} \geq 1 + \frac{d(x,A_i)}{4}$. Now we return our attention to $H'_i(x)-H'_i(y)$, which we deduce from our previous arguments can be written in the following way:
\[
 H'_i(x)-H'_i(y) = \sum \mu_n\psi_i(b_n),
\]
where each $b_n \in \partial_i(x)\cup \partial_i(y)$ and $\mu_n$ is the value of the function $(\ref{diffl1})$ at $b_n$. From the above argument we know that $\sum \mu_n=0$ and $\sum \abs{\mu_n} \leq \frac{C'}{d(x,A_i)+1}$.
\m
But for any two points $a,b\in \partial_i(x)\cup \partial_i(y)$, $\norm{\psi_i(a)-\psi_i(b)}\leq 2K$, by conditions (C1) and (C2) and the fact that each $\psi_i$ is $1$-Lipschitz. Therefore,
\[
 \norm{H'_i(x)-H'_i(y)}_p \leq \frac{2KC'}{d(x,A_i)+1}.
\]
Instead, assume without loss of generality that $k_{x,i}>a_{x,i}$, then $\partial_i(\geo{G_{x,k}})=\emptyset$ for some $k$ and by condition (C2) the length of any $i$-domain of any considered geodesic is bounded from above by $K$.
\m
Hence, using the fact that $\abs{k_{x,i}-k_{y,i}}\leq K$, we deduce in the same way as above that $(\ref{diffl1})$ has $\ell^1$ norm bounded by $\frac{C'}{d(x,A_i)+1}$ for some uniform constant $C'$.
\m
Again writing
\[
 H'_i(x)-H'_i(y) = \sum \mu_n\psi_i(b_n),
\]
we see that as each $\psi_i$ is $1$-Lipschitz,
\[
 \norm{H'_i(x)-H'_i(y)}_p \leq \sum \abs{\mu_n}\norm{\psi_i(b_n)}_p \leq \frac{2KC'}{d(x,A_i)+1},
\]
completing the lemma. \qed

\subsection{The proof of Theorem $\ref{atg'}$}\label{pfatg}
Now we are ready to prove Theorem $\ref{atg'}$ using the embedding
\[
\phi:X\to \ell^p(X)\oplus \bigoplus_{i\in I} \ell^p(X_i) \h \h \tu{given by} \h \h \phi(x)=\phi^s(x)+\phi^l(x).
\]
This is Lipschitz by Lemmas \ref{phi^sLip} and \ref{phi^lLip}.
\m
Consider $x,y\in X$ with $d(x,y)\geq CK$ ($C$ is chosen such that $\rho'(CK)\geq 35K$ and $C\geq 35$).
\m
Fix geodesics $\geo{g_x}\in\setg{x}{e}$ and $\geo{g_y}\in\setg{y}{e}$.
\m
Set $x_y$ to be the closest point $p_{x,y}$ on $\geo{g_x}$ to $e$ such that $d(p_{x,y},\geo{g_y})\geq 5K$ and define $y_x$ similarly. 
\m
Notice that if $x_y,y_x\in A_i$ for some $i$, then that $i$ is unique, as the intersection of any two pieces has diameter at most $K$, by condition (C3).
\m
Let $J_x=\setcon{j\in I_x}{\geo{g_x}|_{\setg{x}{x_y}} \cap A_j \neq \emptyset}$ and $J'_x = J_x\cap I'_x$. We define $J_y$ and $J'_y$ similarly.
\m
$J_x\cap J_y$ has cardinality at most $1$, by condition (C1). 
\m
We now deal with the cases where there $x$ and $y$ are separated by a large piece. To detect this, we will be using the embedding $\phi^l$.
\m
Suppose $\abs{J_x\cap J_y}=1$, label that index $i$ and suppose $d(\geo{g_x}|_i^+, \geo{g_y}|_i^+)\geq \frac{1}{7} d(x,y)$, then
\[
\norm{\phi(x)-\phi(y)}_p^p \geq \norm{H'_i(x)-H'_i(y)}_p^p.
\]
We notice that the sets $\partial_i(x)$ and $\partial_i(y)$ are disjoint as $\partial_i(x)$ has diameter at most $K$, so the function defined in the proof of Lemma \ref{phi^lLip}, $(\ref{diffl1})$ has $\ell^1$ norm $2$ in this case. Therefore, we can write
\[
H'_i(x)-H'_i(y) = \sum_n \mu_n H'_i(x).\chi(\set{b_n}) - \sum_m \mu_m H'_i(y)\chi(\set{b_m})
\]
with $\mu_a$ again being the value of $(\ref{diffl1})$ evaluated at $b_a$. Notice that $\sum_n \mu_n = \sum_m \mu_m = 1$ and the sets $\set{b_m}$ and $\set{b_n}$ are disjoint. Pairing up the $\mu_n$ and $\mu_m$ and applying condition (C2) we see that
\[
\norm{H'_i(x)-H'_i(y)}_p \geq \rho'(d(\geo{g_x|_i^+},\geo{g_y|_i^+})) - 4K \geq \frac{1}{35}\rho'(d(x,y)),
\]
where the last step comes from the concavity of $\rho'$ and the upper bound on $d(x,y)$.
\m
If $J_x\cap J_y = \set{i}$ we now set $x_y=\geo{g_x}|_i^+$ and $y_x=\geo{g_y}|_i^+$. Otherwise we leave $x_y$ and $y_x$ as before. In particular, $d(x_y,y_x)\geq \frac{6}{7}d(x,y) - 2K$.
\m
Without loss of generality, suppose that $d(x,x_y)\geq d(y,y_x)$. From this we deduce that $d(x,x_y)\geq \frac{3}{7}d(x,y)-K$. If there exists some $j\in J_x\setminus J_y$ with $l_j(\geo{g_x})\geq \frac{1}{7}d(x,y)$, then
\[
\begin{array}{rcl} 
\displaystyle \norm{\phi(x)-\phi(y)}_p
&
 \geq 
&
\norm{H'_j(x)-H'_j(y)}_p = \norm{H'_j(x)}_p 
\m
&
\geq 
&
\rho'(l_j(\geo{g_x})) - 2K \geq \frac{1}{35}\rho'(d(x,y)). 
\end{array}
\]
If this does not happen, then we use the embedding $\phi^s$ to detect the distance between $x$ and $y$.
\[
\norm{\phi(x)-\phi(y)}_p^p \geq \sum_{j\in J'_x\setminus J_y} \frac{f(d(x,A_i))^p}{d(x,A_i)} \norm{H_j(x)}_p^p.
\]
As every point $p$ lying on $\geo{g_x}$ at distance between $\frac{2}{3}d(x,x_y)$ and $d(x,x_y)$ from $x$ lies in some $A_j$ with $j\in J'_x\setminus J_y$, we use Lemma \ref{lbound} and the lower bound of Lemma \ref{2'} to deduce
\[
\norm{\phi(x)-\phi(y)}_p^p \succeq f\left(\frac{1}{3}d(x,x_y)-1\right)^p.
\]
Finally, $f$ is concave and $d(x,x_y)\geq \frac{3}{7} d(x,y)-K$, so this gives
\[
\norm{\phi(x)-\phi(y)}_p \succeq f\left(\frac{1}{7}d(x,y) - (K+1)\right) \succeq f(d(x,y)).
\]
\qed

\begin{cor}\label{RHgps'} \h \tu{(cf. Theorem \ref{RHgps})}
\m
Let $G$ be a finitely generated group which is hyperbolic relative to a collection of finitely generated subgroups $\set{H_i}$. Given any $p>1$, any collection of coarse embeddings of $H_i$ into $\ell^p$ spaces with associated concave lower bounding functions $\rho^i_-$ and any function $f$ with property \tu{($C_p^c$)} there is a map $\phi$ from $G$ into an $\ell^p$ space such that for all $x,y\in G$,
\[
 \rho(d(x,y))\preceq \norm{\phi(x)-\phi(y)}_p \preceq d(x,y),
\]
where $\rho(n) = \min\set{\rho^i_-(n), f(n)}.$
\m
In particular, for all $p>1$, $\alpha_p^*(G)=\min\set{\alpha_p^*(H_i)}$.
\end{cor}
\pf This follows from Theorem \ref{atg'}, Proposition \ref{rhinX} and Appendix A of \cite{DS05}.
\m
In particular, we obtain an embedding result for all closed $3$-manifolds.

\begin{cor}\label{3mgps} \h Let $M$ be a closed $3$-manifold, then for all $p>1$ and all $f$ satisfying property \tu{($C_p^c$)}, there exists a map $\phi$ from $\pi_1(M)$ into some $\ell^p$ space, such that for all $x,y\in\pi_1(M)$,
\[
 f(d(x,y))\preceq \norm{\phi(x)-\phi(y)}_p\preceq d(x,y).
\]
So $\alpha_p^*(\pi_1(M))=1$ for all $p>1$.
\end{cor}
\pf Consider first the geometric manifolds. The fundamental groups of these are quasi-isometric to one of eight Thurston geometries, which are either compact, Euclidean, hyperbolic or have a suitable embedding via \cite[Theorem $1$]{Te11}.
\m
In the non-geometric case, we decompose the manifold along tori using the Geometrisation Theorem \cite{Per1,Per2,CZ06a,CZ06b,KL08,MT07,MT08}. If $M$ has no hyperbolic part then it is a graph manifold and Smirnov proves this has finite Assouad-Nagata dimension, \cite{Sm10}. A theorem of Gal then gives suitable bounds on compression \cite{Gal}. For an explicit embedding, one can use Theorem \ref{TG'} and \cite[Theorem $1$]{HS11}.
\m
Finally, if it has a hyperbolic part, then $\pi_1(M)$ is hyperbolic relative to the fundamental groups of a finite collection of graph manifold groups and virtually polycyclic groups \cite{Da03}. Using \cite[Theorem $1$]{Te11} and \cite{Gal} again and applying Corollary $\ref{RHgps}$ completes the result. \qed
\m
Since this paper was originally written, this result has been improved by Mackay-Sisto \cite{MS12}, who prove that closed $3$-manifold groups have finite Assouad-Nagata dimension, so the above result follows from \cite{Gal}.
\m
By way of complete contrast, Sapir proves the existence of a closed aspherical $4$-manifold $M$ where $\pi_1(M)$ coarsely contains expanders and hence admits no coarse embedding into any Hilbert space \cite{Sa11}. This uses Gromov's proof \cite{Gr00} of the existence of a finitely generated group coarsely containing expanders.
\m
Finally, we also obtain an estimate for $L_p$ compression.

\begin{cor}\label{L_p} \h Let $X$ be the $0$-skeleton of an asymptotically tree-graded simplicial graph of uniformly bounded degree and let $\setcon{A_i}{i\in I}$ be a suitable choice of pieces. Suppose we are given a collection of maps $\psi_i:A_i\to L_p([0,1])$ and a concave function $\rho':\N\to\R_{\geq 0}$ such that for all $i\in I$ and all $x,y\in A_i$,
\[
 \rho'(d(x,y))\leq \norm{\psi_i(x)-\psi_i(y)}_p \leq d(x,y).
\]
For each function $f$ satisfying property \tu{($C_q^c$)} where $q=\max\set{p,2}$ there exists a map $\phi$ of $X$ into $L_p([0,1])$ with
\[
\min\{f(d(x,y)),\rho'(d(x,y))\}\preceq \norm{\phi(x)-\phi(y)} \preceq d(x,y).
\]
\end{cor}
\pf $X$ is countable, and we enumerate it as $X=\set{x_1,x_2,\dots}$. The map $\phi:\ell^p(X) \into L^p([0,1])$ defined as the linear extension of the mapping $x_n \mapsto f_n$ given by
\[
 f_n(x) = \left\{
 \begin{array}{rl}
  2^\frac{n+1}{p}
  &
  \tu{if } x\in (2^{-(n+1)},2^{-n})
  \m
  0
  &
  \tu{otherwise}
 \end{array} \right.
\] 
is an isometric embedding of $\ell^p(X)$ into $L^p([0,1])$. The remainder then follows by recalling the fact that $L^2([0,1])$ isometrically embeds into $L_p([0,1])$ when $p\in[1,2]$, \cite{Woj} and applying Theorem \ref{atg'}. \qed

{
\footnotesize
\newcommand{\etalchar}[1]{$^{#1}$}

}
\end{document}